\documentclass{dcds-s}
\usepackage{psfrag}
\usepackage{amsmath}
  \usepackage{paralist}
  \usepackage{graphics} 
  \usepackage{epsfig} 
\usepackage{graphicx}  \usepackage{epstopdf}
 \usepackage[colorlinks=true]{hyperref}
\hypersetup{urlcolor=blue, citecolor=red}

\def\currentvolume{}
 \def\currentissue{}
  \def\currentyear{}
   \def\currentmonth{}
    \def\ppages{}
     \def\DOI{}

\newtheorem{theorem}{Theorem}[section]

\theoremstyle{definition}

\newcommand{\R}{\mathbb R}
\newcommand{\Dx}{\Delta x}

\newcommand{\rhomax}{\rho_{\textup{max}}}
\newcommand{\DtsuDx}{\frac{\Delta t}{\Delta x}}

\title[An easy algorithm for traffic flow on networks]
      {An easy-to-use algorithm for simulating traffic flow on networks: Numerical experiments}

\author[Gabriella Bretti, Maya Briani and Emiliano Cristiani]{}


 \email{g.bretti@iac.cnr.it}
 \email{m.briani@iac.cnr.it}
 \email{e.cristiani@iac.cnr.it}

\thanks{The research leading to these results has received funding from the
European Union FP7 under grant No.\ 257462 HYCON2 Network of
Excellence. The research was also supported by the collabo\-ration with the private company ZEROPIU and by the Google Research Award
``Multipopulation Models for Vehicular Traffic and Pedestrians'', 2012-2013.}

\begin{document}
\maketitle

\centerline{\scshape Gabriella Bretti, Maya Briani and Emiliano Cristiani}
\medskip
{\footnotesize
 \centerline{Istituto per le Applicazioni del Calcolo ``M. Picone''}
   \centerline{Consiglio Nazionale delle Ricerche}
   \centerline{Via dei Taurini, 19 -- 00185 Rome, Italy}
} 



%

\begin{abstract}
In this paper we propose a Godunov-based discretization of a hyperbolic system of conservation laws with discontinuous flux, modeling vehicular flow on a network. Each equation describes the density evolution of vehicles having a common path along the network.
We show that the algorithm selects automatically an admissible solution at junctions, hence \emph{ad hoc} external procedures (e.g., maximization of the flux via a linear programming method) usually employed in classical approaches are no needed. Since users have not to deal explicitly with vehicle dynamics at junction, the numerical code can be implemented in minutes.
We perform a detailed numerical comparison with a Godunov-based scheme coming from the classical theory of traffic flow on networks which maximizes the flux at junctions.
\end{abstract}

\section{Introduction}\label{sec:introduction}
Starting from the introduction of the LWR model \cite{LW, R}, a huge literature about macroscopic fluid-dynamic models for traffic flow was developed. More recently, models, theory and numerical approximations for traffic flow on networks became a hot topic \cite{coclite2005SIMA, daganzo1995TRBa, piccolibook, holden1995SIMA}. The interest in forecasting traffic flow on large networks became even stronger in the very last years, due to the increasing number of GPS devices (smartphones, satellite navigators, black boxes) which provide real-time traffic data. Private companies like GOOGLE, WAZE MOBILE, NOKIA \cite{herrera2010TRB}, INRIX, OCTOTELEMATICS \cite{cristiani2010CAIM}, ZEROPIU, YANDEX, started collecting data and, in some cases, broadcasting traffic forecast.

\emph{The model.} In this paper we study from the numerical point of view a first-order version of the model proposed in \cite[Sect.4]{hilliges1995TRB}. Let us consider a network, i.e.\ a directed graph with $N_R$ arcs (roads) and $N_J$ nodes (junctions). \textit{Vehicles moving on the network are divided on the basis of their desired path}. Let us assume that the number of possible paths on the graph is $N_P$ and denote those paths by $P^1,\ldots,P^p,\ldots,P^{N_P}$. We stress that paths can share some arcs of the networks. A point $x^{(p)}$ of the network is characterized by both the path $p$ it belongs to and the distance $x$ from the origin of that path. We denote by $\mu^p(x^{(q)},t)$ the density of the vehicles following the $p$-th path at point $x^{(q)}$ at time $t>0$, and we assume that $\mu^p(x^{(q)},t) \in [0,\rhomax]$ for some maximal density $\rhomax$. Note that we have, by definition, $\mu^p(x^{(q)},t)=0$ if $x^{(q)}\notin P^p$.
We also define
\begin{equation}\label{def:omega}
\omega(x^{(p)},t):=\sum_{q=1}^{N_P}\mu^q(x^{(p)},t),
\end{equation}
i.e.\ $\omega(x^{(p)},t)$ is the sum of all densities living at $x^{(p)}$ at time $t$.
The function $\omega$ is discontinuous and takes into account the topology of the network.
Note that, for any point $x^{(q)}$, the densities $\mu^p(x^{(q)},t)$, $p=1,\ldots,N_P$, are admissible if their sum $\omega(x^{(q)},t)\leq\rhomax$.
Let us denote by $v(\omega)$ the velocity of vehicles (given as a function of the density) and by $f(\omega)=\omega v(\omega)$ the flux of vehicles.
The LWR-based model is constituted by the following system of $N_P$ conservation laws with space-dependent and discontinuous flux
\begin{equation}\label{eq:system_v}
\frac{\partial}{\partial t}\mu^p(x^{(p)},t)+
\frac{\partial}{\partial x^{(p)}}\left(\mu^p(x^{(p)},t)\ v\big(\omega(x^{(p)},t)\big)\right)=0, \quad x^{(p)}\in P^p,\ t>0,
\end{equation}
or, equivalently,
\begin{equation}\label{eq:system_f}
\frac{\partial}{\partial t}\mu^p(x^{(p)},t)+
\frac{\partial}{\partial x^{(p)}}\left(\frac{\mu^p(x^{(p)},t)}{\omega(x^{(p)},t)}\ f\big(\omega(x^{(p)},t)\big)\right)=0, \quad x^{(p)}\in P^p,\ t>0,
\end{equation}
for $p=1,\ldots,N_P$.
If $\omega=0$ we have, \textit{a fortiori}, $\mu^p=0$, then it is convenient to set $\frac{\mu^p}{\omega}=0$ in (\ref{eq:system_f}) to avoid singularities.
In the following we assume that the flux $f\in C^0([0,\rhomax])\cap C^1((0,\rhomax))$ and
\begin{equation}\label{proprieta_f}
f(0)=f(\rhomax)=0,\qquad f \textrm{ is strictly concave}, \qquad f(\sigma)=\max_{\omega\in (0,\rhomax)}f(\omega).
\end{equation}
Equations of the system (\ref{eq:system_v}) (or (\ref{eq:system_f})) are coupled by means of the velocity $v$, which depends on the total density $\omega$. On the other hand, not all the equations of the system are coupled with each other because paths do not have necessarily arcs in common.

It is plain that the number of equations in the system grows rapidly when the number of nodes of the graph increases, making unfeasible the computation of a numerical solution. To keep the computational load within reasonable limits, we also propose another version of the model which splits the vehicles on the basis on their path only at junctions. In this way, along the arcs a single equation for the total density $\omega$ is considered. The price to pay is that the global behavior of drivers is lost.


\emph{Relevant literature.} The \textit{multi-path} model described above differs from the so-called \emph{multi-population} or \emph{multi-class} models, see, e.g., \cite{benzoni2003EJAM, wong2002TRA}. In those cases, the models consist of one equation for a single road (extension to networks is also possible) with different velocity functions $v_i$, one for each class of vehicles. Typically, the populations have different maximal velocities, in order to take into account different types of vehicles or drivers' behaviors.

An apparently similar model is the one presented in \cite{garavello2005CMS} (see also \cite{lebacque1996proc}). In that case, vehicles are divided in different populations on the basis on their source-destination pair. Given the total density $\omega$ of all vehicles, the density $\mu^p$ of the $p$-th population is given by
$$
\mu^p(x,t)=\pi^p\big(x,t,O(p),D(p)\big)\ \omega(x,t),
$$
where $\pi^p(x,t,O,D)$ specifies the percentage of the total density that started from source $O$, it is moving towards the final destination $D$, and it is at $x$ at time $t$. Moreover, $O(p)$, $D(p)$ are the origin and the destination associated to the $p$-th path, respectively. Then, a standard PDE for $\omega$ is coupled with a system of PDEs for the coefficients $\pi^p$'s.

Several papers investigate from the theoretical point of view the (systems of) scalar conservation laws. The interested reader can find in the book \cite{levequebook} an introduction to the field, in \cite{adreianov2011ARMA,burger2008JEM} some references for the case of discontinuous flux, and in the book \cite{bressanbook} the analysis of the systems of conservation laws. Systems of scalar conservation laws with discontinuous flux are instead less studied. An attempt related to traffic flow can be found in \cite{mercier2009JMAA}, where a model very similar to the one considered here is investigated.
From the numerical point of view, a good basic reference is again the book \cite{levequebook}. We also point out the paper \cite{herty2007ANM}, where a numerical method for (systems of) scalar conservation laws with discontinuous flux is proposed, and the paper \cite{towers2000SINUM}, where the convergence of a Godunov-based scheme for scalar conservation laws with discontinuous flux is investigated.


\emph{Goal.} The goal of this paper is proposing an algorithm to solve numerically the system (\ref{eq:system_f}). The numerical method is based on the Godunov scheme. The algorithm does not require use of \emph{ad hoc} procedures (e.g., linear programming) to solve the dynamics at junctions, as it is done in classical approaches \cite{bretti2007ACME, daganzo1995TRBa, piccolibook, herty2007ANM}. This leads to a very simple algorithm which can be implemented in minutes. We perform a detailed comparison with the algorithm described in \cite{bretti2007ACME, piccolibook} based on the LWR model and the Godunov scheme. We show that the two algorithms do not coincide in general since the proposed one does not always maximize the flux at junction. Nevertheless, the solution is admissible in the framework of the classical theory and in some cases it is expected to provide a more reasonable behavior of the car flow. A theoretical investigation of the algorithm proposed here can be found in \cite{briani2013inprep}.


\emph{Paper organization.} The paper is organized as follows:
In Section 2 we recall some basic facts of the classical theory of traffic flow on networks and Godunov discretization, following \cite{bretti2007ACME, coclite2005SIMA, piccolibook}.
In Section 3 we present the algorithm for small and large networks.
In Sections 4, 5 and 6 we report the numerical results and comparisons for the junctions with 2 incoming roads and 1 outgoing road, 1 incoming road and 2 outgoing roads, and 2 incoming roads and 2 outgoing roads, respectively.
In Section 7 we report the result of an experiment performed on a real network in Rome, Italy.
Finally, in Section 8 we sketch some conclusions.


\section{Classical theory}\label{sec:classicaltheory}
In this section we recall some basic facts of the classical theory of traffic flow on networks, including a numerical approximation based on the Godunov scheme. Our main references are \cite{bretti2007ACME, coclite2005SIMA, piccolibook}.

\subsection{Basic definitions, assumptions, and results}
Let $I=[a,b]$ be a generic arc of the graph, i.e.\ a road.
At any time $t$, the evolution of the vehicle density on the network is computed by a two-step procedure: First, a classical conservation law is solved at any internal point of the arcs; Second, the densities at endpoints $a$, $b$ which correspond to a junction are computed. The latter step has not in general a unique admissible solution, so that additional constraints must be added. Beside the conservation of vehicles at junctions, we assume here that drivers behave in order to maximize the flux at junctions and that incoming roads are regulated by priorities (right of way). The second step is performed by a linear programming method.

On each arc, the density $\rho(x,t)$ of all vehicles (no distinction among vehicles is considered here) is given by the entropic solution of
$$
\frac{\partial}{\partial t}\rho+\frac{\partial}{\partial x}f(\rho)=0, \qquad x\in I, \quad t>0.
$$
We consider now a generic junction $J$, with $I_1,\ldots,I_n$ incoming roads and $I_{n+1}$, $\ldots$, $I_{n+m}$ outgoing roads. We assume that the choice of the outgoing road is prescribed by a \textit{matrix of preferences},
\begin{equation}\label{def:A}
A=
\left(
\begin{array}{ccc}
\alpha_{n+1,1} & \ldots & \alpha_{n+1,n} \\
\ldots & \ldots & \ldots \\
\alpha_{n+m,1} & \ldots & \alpha_{n+m,n}
\end{array}
\right)
\end{equation}
where $0\leq \alpha_{j,i}\leq 1$ for every $i\in\{1,\ldots,n\}$ and $j\in\{n+1,\ldots,n+m\}$, and
$\sum_{j=n+1}^{n+m}\alpha_{j,i}=1$, $i=1,\ldots,n$.
The $i$-th column of $A$ describes how the traffic from $I_i$ distributes in percentage to the outgoing roads.

A basic requirement for a vector $(\rho^1,\ldots,\rho^{n+m})$ to be an admissible solution to the problem at $J$ is that
\begin{equation}\label{RH}
\sum_{i=1}^n f(\rho^i(b_i-,t))=\sum_{j=n+1}^{n+m} f(\rho^j(a_j+,t)),
\end{equation}
which translates the fact that the vehicles are conserved at junction. Note that (\ref{RH})
can be seen as a generalization of the Rankine-Hugoniot condition at junctions.

We define
\begin{equation}\label{def:gammamaxin}
\gamma^i_{\text{max}}(\rho(b_i,t)) =
\left\{\begin{array}{ll}
f(\rho(b_i,t)), & \mbox{ if } \rho(b_i,t)\in[0,\sigma],\\
f(\sigma), & \mbox{ if } \rho(b_i,t)\in]\sigma,\rhomax],
\end{array}\right. \quad i=1,\ldots,n
\end{equation}
and
\begin{equation}\label{def:gammamaxout}
\gamma^j_{\text{max}}(\rho(a_j,t)) =
\left\{\begin{array}{ll}
f(\sigma), & \mbox{ if } \rho(a_j,t)\in[0,\sigma],\\
f(\rho(a_j,t)), & \mbox{ if } \rho(a_j,t)\in]\sigma,\rhomax],
\end{array}\right. \quad j=n+1,\ldots,n+m.
\end{equation}
The quantities $\gamma^i_{\textup{max}}(\rho(b_i,t))$ and $\gamma^j_{\textup{max}}(\rho(a_j,t))$ represent, respectively, the maximum incoming and the maximum outgoing flux that can be obtained on each road. Then we define
\begin{eqnarray}
& &\hskip-30pt \Omega_i:= [0,\gamma^i_{\text{max}}(\rho(b_i,t))],\qquad i=1,\ldots, n, \label{Omegain} \\
& &\hskip-30pt \Omega_j:= [0,\gamma^j_{\text{max}}(\rho(a_j,t))],\qquad j=n+1,\ldots, n+m, \label{Omegaout} \\
& &\hskip-30pt \Omega:=\{(\gamma^1,\ldots,\gamma^n)\in\Omega_1\times\ldots\times\Omega_n\ | \ A (\gamma^1,\ldots,\gamma^n)^T\in\Omega_{n+1}\times\ldots\times\Omega_{n+m}\} \label{Omega}.
\end{eqnarray}
The sets $\Omega_i$, $\Omega_j$ contain all the possible fluxes for the solution at junctions and then the set $\Omega$ contains all the possible admissible fluxes at the end of the incoming roads, taking into account the matrix of preferences $A$.
Since we want to maximize the flux at junction, we find the solution(s) of the maximization problem with linear constraints
\begin{equation}\label{PL}
\max_{(\gamma^1,\ldots,\gamma^n)\in\Omega}\sum_{i=1}^n \gamma^i.
\end{equation}
Note that the problem (\ref{PL}) has not in general a unique solution. For example, if $n=2$ and $m=1$ (two incoming and one outgoing roads) we have $A=(1 \ 1)$, and the constraint reads as $\gamma^1+\gamma^2\leq\gamma^{3}_{\text{max}}$. If $\gamma^1_{\text{max}}+\gamma^2_{\text{max}}>\gamma^3_{\text{max}}$, we have infinite solutions.
To fix this, the additional constraint
\begin{equation}\label{constraint_priorities}
(\gamma^1,\ldots,\gamma^n)\in\{qs,\ s\in\R^+\},
\end{equation}
can be introduced, where $q=(q_1,\ldots,q_n)$ are the so-called \textit{priorities} (right of way) coefficients such that $q_i\geq 0\ \forall i$ and $\sum_i q_i=1$.
Equation (\ref{constraint_priorities}) translates the fact that some incoming roads have priority with respect to the other roads in assigning their flux to the outgoing roads. The constraint \eqref{constraint_priorities} guarantees uniqueness of the solution of problem (\ref{PL}) (unless a further projection onto $\Omega$ is needed, see \cite{chitour2005DCDS-B,piccolibook} for details).
Finally, define $(\gamma^*_1,\ldots,\gamma^*_n)$ as the unique solution of problem (\ref{PL}) (with the additional constraint  (\ref{constraint_priorities}) if necessary),
\begin{equation}\label{def:gammastar_i}
(\gamma^*_1,\ldots,\gamma^*_n):= \operatorname*{arg\,max}_{(\gamma^1,\ldots,\gamma^n)\in\Omega} \sum_{i=1}^n \gamma^i
\end{equation}
and, consequently,
\begin{equation}\label{def:gammastar_j}
\gamma^*_{j}:=\sum_{i=1}^{n}\alpha_{j,i} \ \gamma^*_i,\qquad j=n+1,\ldots, n+m.
\end{equation}


\subsection{Numerical approximation by the Godunov scheme}
Let us describe briefly the Godunov scheme for solving a conservation law of the form
$$
\left\{
\begin{array}{l}
\displaystyle{\frac{\partial}{\partial t}\rho+\frac{\partial}{\partial x}f(\rho)=0}\\[3mm]
\rho(x,0)=\bar\rho(x)
\end{array}
\right.
$$
on a generic arc $I=[a,b]\subset\R$ of the network and $t\in[0,t_f]$, where $t_f$ is the final time. We define a numerical grid in $[a,b]\times[0,t_f]$ with space step $\Delta x$ and time step $\Delta t$, satisfying the CFL condition
\begin{equation}\label{CFLmodelloclassico}
\Delta t \max_{\rho}|f'(\rho)|\leq \Delta x.
\end{equation}
Let us denote by $(x_k$,\ $t^n):=(k\Delta x,\ n\Delta t )$, $k\in\mathbb Z$, $n\in\mathbb N$, the generic grid node and by $\rho_k^n$ the density $\rho$ at $(x_k,t^n)$.
Once the initial condition $\bar\rho$ has been projected on the grid,
in the internal nodes of the interval $[a,b]$ the density at time $t^{n+1}$ is given by
\begin{equation}\label{Godunov}
\rho_k^{n+1}=\rho_k^n-\frac{\Delta t}{\Delta x}\Big(g(\rho_k^n,\rho_{k+1}^n)-g(\rho_{k-1}^n,\rho_{k}^n)\Big),
\end{equation}
where the numerical flux $g$ is defined as
\begin{equation}\label{GodunovFlux}
g(\rho_-,\rho_+):=
\left\{
\begin{array}{ll}
\min\{f(\rho_-),f(\rho_+)\} & \textrm{if } \rho_-\leq \rho_+ \\
f(\rho_-) & \textrm{if } \rho_->\rho_+ \textrm{ and } \rho_-<\sigma \\
f(\sigma) & \textrm{if } \rho_->\rho_+ \textrm{ and } \rho_- \geq \sigma \geq \rho_+ \\
f(\rho_+) & \textrm{if } \rho_->\rho_+ \textrm{ and } \rho_+>\sigma
\end{array}
\right. .
\end{equation}
At the boundary nodes we proceed as follows:
\begin{itemize}
\item If the node is not linked to any junction, then we assign the desired boundary condition (Dirichlet or Neumann).
\item If the node is a right endpoint and corresponds to the $i$-th incoming road of a junction, we use the maximal flux (\ref{def:gammastar_i}) and set
\begin{equation}\label{godunov_ben_out}
\rho_k^{n+1}=\rho_k^n-\frac{\Delta t}{\Delta x}\Big(\gamma^*_i-g(\rho_{k-1}^n,\rho_{k}^n)\Big).
\end{equation}
\item If the node is a left endpoint and corresponds to the $j$-th outgoing road of a junction, we use the maximal flux (\ref{def:gammastar_j}) and set
\begin{equation}\label{godunov_ben_in}
\rho_k^{n+1}=\rho_k^n-\frac{\Delta t}{\Delta x}\Big(g(\rho_k^n,\rho_{k+1}^n)-\gamma^*_j\Big).
\end{equation}
\end{itemize}


\section{Numerical approximation of the multi-path model on small and large networks}\label{sec:numerical approximation}
In this section we present a numerical approximation for the system (\ref{def:omega}),(\ref{eq:system_f}) based on the Godunov scheme (\ref{Godunov}). Then, we discuss how the multi-path model can be modified to deal easily with large networks.

\subsection{Numerical approximation}\label{sec:num_approx}
Let us denote by $\mu^{n,p}_{k^{(q)}}$ the approximate density $\mu^p(x_k^{(q)},t^n)$, where $k^{(q)}$ is the $k$-th node along the path $P^q$. Then, analogously to (\ref{def:omega}), we define
\begin{equation}\label{def:omega_approx}
\omega^{n}_{k^{(p)}}:=\sum_{q=1}^{N_P}\mu^{n,q}_{k^{(p)}}, \qquad q=1,\ldots,N_P.
\end{equation}
From now on, to avoid cumbersome notations, we write $k^p$ instead of $k^{(p)}$.
Computing properly $\omega^{n}_{k^p}$ at every node is the sole part of the algorithm which requires some effort. Indeed, the rest of the algorithm is constituted by the imposition of the boundary conditions at the beginning and the end of the $N_P$ paths and by the computation of the discrete solution at any internal nodes $k^p$ by means of the following conservative scheme:
\begin{equation}\label{schema}
\mu_{k^p}^{n+1,p} = \mu_{k^p}^{n,p}-
\DtsuDx\left(\frac{\mu_{k^p}^{n,p}}{\omega_{k^p}^{n}} \ g(\omega_{k^p}^{n},\omega_{k^p+1}^{n})-
\frac{\mu_{k^p-1}^{n,p}}{\omega_{k^p-1}^{n}}\ g(\omega_{k^p-1}^{n},\omega_{k^p}^{n})
\right)
\end{equation}
for $n>0$ and $p=1,\ldots,N_P$.
Note the intrinsic asymmetry of (\ref{schema}). The coefficients in front of the fluxes involve only the nodes $k^p$ and $k^p-1$, and not $k^p+1$.

We stress again that no special management of the junctions is needed. \textit{The simplicity of the scheme is the main strength of this approach, making it possible to simulate traffic flow on networks in minutes}.
%
%
\subsection{A local model for large networks}\label{sec:localmodel}
If the network under consideration is small, the number of possible paths is limited. Then, the number of equations in the system (\ref{def:omega}),(\ref{eq:system_f}) fits a manageable size. This is true even if we deal with large networks where many paths are impracticable or negligible. Conversely, if the network is large and has a large number of paths, the computation of $\omega^{n}_{k^p}$ becomes a hard task. In this case, we adopt a hybrid point of view, creating an algorithm which merges the features of the multi-path algorithm described above and the classical algorithm described in Section \ref{sec:classicaltheory}. In the internal nodes of the arcs a single equation for the \emph{total} density $\rho$ (or $\omega$) is solved. Then, just before each junction, vehicles are split on the basis of their direction by means of the matrix of preferences $A$, see (\ref{def:A}), and the scheme \eqref{schema} is applied. Just after the junction, densities associated to the same arc are summed again. Considering that typical real junctions have at most three incoming and three outgoing roads, we have at most nine different paths at each junction.

\begin{figure}[h!]
\begin{center}
\begin{psfrags}
\psfrag{J1}{$J_1$} \psfrag{J2}{$J_2$}
\psfrag{I1}{$[a_1,b_1]$} \psfrag{I2}{$[a_2,b_2]$} \psfrag{I3}{$[a_3,b_3]$} \psfrag{I4}{$[a_4,b_4]$} \psfrag{I5}{$[a_5,b_5]$}
\includegraphics[width=0.7\textwidth]{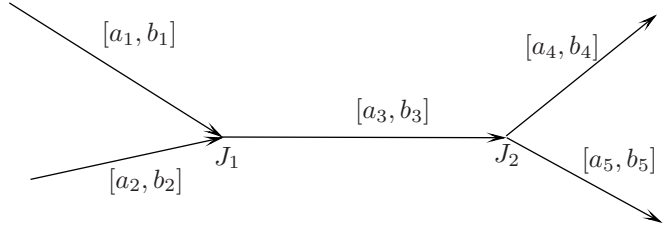}
\end{psfrags}
\end{center}
\caption{A network with 5 arcs and 2 junctions.}
\label{fig:esempio_network_nonlocale}
\end{figure}

As already noted in \cite{daganzo1994TRB, garavello2005CMS}, the original (global) multi-path approach of Section \ref{sec:num_approx} \emph{differs} from the hybrid (local) method, because splitting the density just before the junctions causes the lost of the global behavior of the drivers. To better understand the difference, let us consider the network in Fig.\ \ref{fig:esempio_network_nonlocale}.
Let us consider only two paths,
$$P^1=[a_1,b_1]\cup [a_3,b_3]\cup [a_4,b_4] \textrm{ and } P^2=[a_2,b_2]\cup [a_3,b_3]\cup [a_5,b_5].$$
If, at some time $t$, vehicles along $[a_2,b_2]$ are stopped by, say, an accident, the two algorithms will result in two different outcomes in $[a_5,b_5]$. The global algorithm will empty the arc $[a_5,b_5]$, while the local algorithm will fill $[a_5,b_5]$ with a percentage of the density in $[a_3,b_3]$. Clearly, only the global algorithm gives the correct solution.
We also note that the classical algorithm described in Section \ref{sec:classicaltheory} is local, while the one presented in \cite{garavello2005CMS} is global.


\section{Two incoming roads and one outgoing road}\label{sec:2in1out}
\begin{figure}[h!]
\begin{center}
\begin{psfrags}
\psfrag{P1}{$P^1$} \psfrag{P2}{$P^2$}
\psfrag{I1}{$[a_1,b_1]$} \psfrag{I2}{$[a_2,b_2]$} \psfrag{I3}{$[a_3,b_3]$}
\psfrag{J}{\tiny $J$} \psfrag{Jp}{\tiny $J\!\!+\!\!1$} \psfrag{Jm}{\tiny $J$-1}
\includegraphics[width=0.45\textwidth]{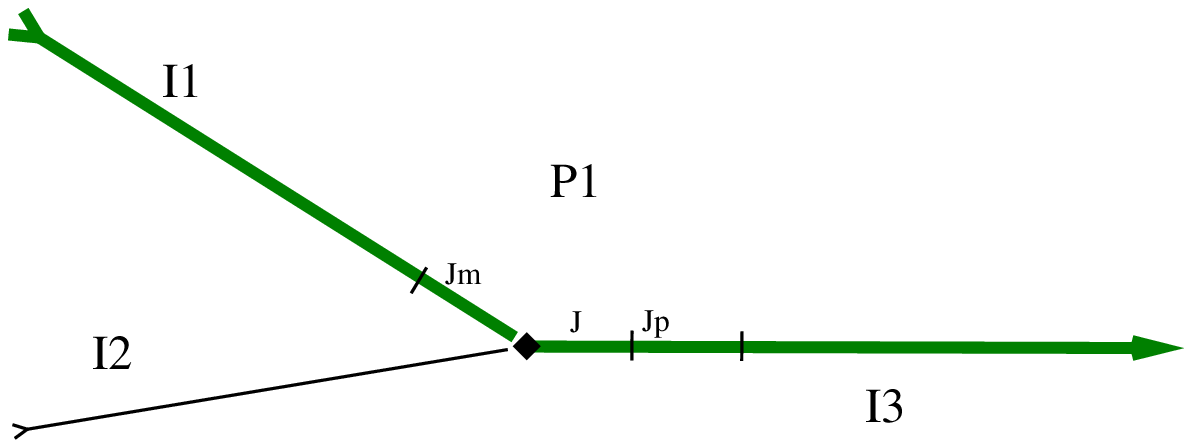}\quad
\includegraphics[width=0.45\textwidth]{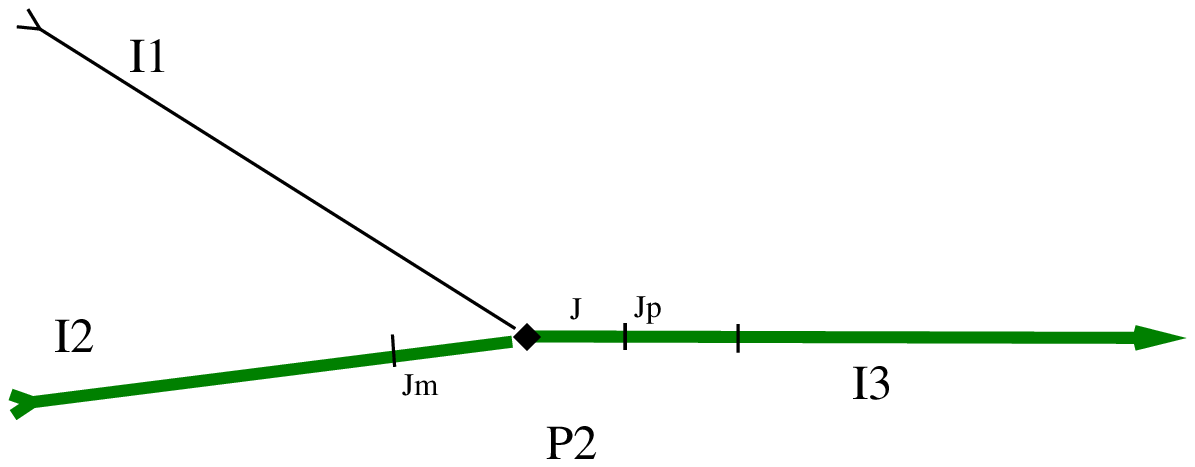}
\end{psfrags}
\end{center}
\caption{A network with 3 arcs and 1 junction. Path $P^1$ (left) and path $P^2$ (right).}
\label{fig:network_2in_1out}
\end{figure}

In this section we consider a network with three arcs and one junction, with two incoming roads and one outgoing road. On this network two paths $P^1$ and $P^2$ are defined, see Fig.\ \ref{fig:network_2in_1out}.
We denote by $J$ the node \emph{just after} the junction, see Fig.\ \ref{fig:network_2in_1out}. Note that the nodes before the junction, namely $J-1$, $J-2$, etc., can refer to one path or the other one, depending on the context.
We have
$$
\omega_{k^1}^{n}=\left\{\begin{array}{ll}\mu_{k^1}^{n,1} & k^1<J \\ \mu_{k^1}^{n,1}+\mu_{k^1}^{n,2} & k^1\geq J\end{array}\right.
,\qquad\qquad
\omega_{k^2}^{n}=\left\{\begin{array}{ll}\mu_{k^2}^{n,2} & k^2<J \\ \mu_{k^2}^{n,1}+\mu_{k^2}^{n,2} & k^2\geq J\end{array}\right.
$$
and the scheme (\ref{schema}) becomes
\begin{equation}\label{schema_2in1out}
\left\{\begin{array}{l}
\mu_{k^1}^{n+1,1} = \mu_{k^1}^{n,1}-\DtsuDx\left(\frac{\mu_{k^1}^{n,1}}{\omega_{k^1}^{n}}g(\omega_{k^1}^{n},\omega_{k^1+1}^{n})- \frac{\mu_{k^1-1}^{n,1}}{\omega_{k^1-1}^{n}}g(\omega_{k^1-1}^{n},\omega_{k^1}^{n})\right),\\
\\
\mu_{k^2}^{n+1,2} = \mu_{k^2}^{n,2}-\DtsuDx\left(\frac{\mu_{k^2}^{n,2}}{\omega_{k^2}^{n}}g(\omega_{k^2}^{n},\omega_{k^2+1}^{n})- \frac{\mu_{k^2-1}^{n,2}}{\omega_{k^2-1}^{n}}g(\omega_{k^2-1}^{n},\omega_{k^2}^{n})\right).
\end{array}\right.
\end{equation}

Let us first prove that the density $\omega_{k^p}^{n}$ is admissible for any $n$ and $k^p$, $p=1,2$, i.e.\ it does not exceed $\rhomax$. This is a crucial property which the scheme must satisfy. It is plain that the risk of exceeding the maximal density is only at node $J$. In the other nodes the result comes from the usual properties of the Godunov scheme.
\begin{theorem}\label{teo:rhoammissibile}
Let the initial densities around the junction be admissible, namely
\begin{equation*}
\mu_{J-1}^{0,1}\leq\rhomax,\quad
\mu_{J-1}^{0,2}\leq\rhomax, \quad
(\mu_{J}^{0,1}+\mu_{J}^{0,2})\leq\rhomax,\quad
(\mu_{J+1}^{0,1}+\mu_{J+1}^{0,2})\leq\rhomax.
\end{equation*}
If the following CFL-like condition holds
\begin{equation}\label{CFLcon2}
2\DtsuDx \sup_{\rho\in(0,\rhomax)}|f'(\rho)|\leq 1,
\end{equation}
then
$$
(\mu_{J}^{n,1}+\mu_{J}^{n,2})\leq\rhomax \qquad \forall n.
$$
\end{theorem}

\noindent\emph{Proof.}
Let us first prove that \eqref{CFLcon2} implies
\begin{equation}\label{condizioneMaya}
\DtsuDx \leq \inf_{\rho\in[\sigma,\! \ \rhomax)}\frac{\rhomax-\rho}{2f(\rho)}.
\end{equation}
Let us define $M:=\sup_{\rho\in(0,\rhomax)}|f'(\rho)|$. By \eqref{CFLcon2} we have
\begin{equation}\label{lambda<=1su2M}
\DtsuDx\leq\frac{1}{2M}.
\end{equation}
Noting that $f(\rhomax)=0$, and using the Lagrange theorem and \eqref{lambda<=1su2M}, we have
$$
\inf_{\rho\in[\sigma,\! \ \rhomax)}\frac{\rhomax-\rho}{2f(\rho)}=
\inf_{\rho\in[\sigma,\! \ \rhomax)}\frac{|\rhomax-\rho|}{2|f(\rho)-f(\rhomax)|}\geq
\inf_{\rho\in[\sigma,\! \ \rhomax)}\frac{1}{2M}=\frac{1}{2M}\geq
\DtsuDx.
$$
To simplify the notations, let us introduce the auxiliary variable
$z^n:=\mu^{n,1}_{J}+\mu^{n,2}_{J}$.
The worst case happens when $\mu_{J-1}^{n,1}=\mu_{J-1}^{n,2}=\sigma$ (incoming roads try to transfer the maximal flux to cell $J$) and $\mu_{J+1}^{n,1}+\mu_{J+1}^{n,2}=\rhomax$ (no flux from cell $J$ to cell $J+1$). In this case the equation for $z$ is
\begin{equation}\label{eq_z}
z^{n+1}=
z^n-\DtsuDx\big(0-g(\sigma,z^n)-g(\sigma,z^n)\big)=
z^n+2\DtsuDx g(\sigma,z^n).
\end{equation}
We proceed by induction: Assume that $z^n\leq \rhomax$ and prove that $z^{n+1}\leq \rhomax$.
We have
$$
g(\sigma,z^n)=
\left\{
\begin{array}{ll}
f(\sigma) & \textrm{ if } z^n\leq\sigma \\
f(z^n)    & \textrm{ if } z^n>\sigma
\end{array}
\right . .
$$
\begin{itemize}
\item CASE 1: $z^n\leq\sigma$\\
We have
$$
z^{n+1}=
z^n+2\DtsuDx f(\sigma)\leq
\sigma+2\DtsuDx f(\sigma).
$$
The conclusion follows easily by \eqref{condizioneMaya}, in particular by the fact that
$$
\DtsuDx \leq \frac{\rhomax-\sigma}{2 f(\sigma)}.
$$
\item CASE 2: $z^n>\sigma$
\begin{itemize}
\item CASE 2.1: $z^n=\rhomax$\\
We have $f(z^n)=f(\rhomax)=0$ and then $z^{n+1}=z^n=\rhomax$.
\item CASE 2.2: $z^n<\rhomax$\\
We have
$$
z^{n+1}=
z^n+2\DtsuDx f(z^n).
$$
The conclusion follows easily by \eqref{condizioneMaya}, in particular by the fact that
$$
\DtsuDx \leq \inf_{\rho\in(\sigma,\! \ \rhomax)}\frac{\rhomax-\rho}{2f(\rho)}.
$$
\end{itemize}
\end{itemize}
The two cases conclude the proof. \hfill $\square$

\medskip

\begin{figure}[h!]
\begin{center}
\includegraphics[width=0.6\textwidth]{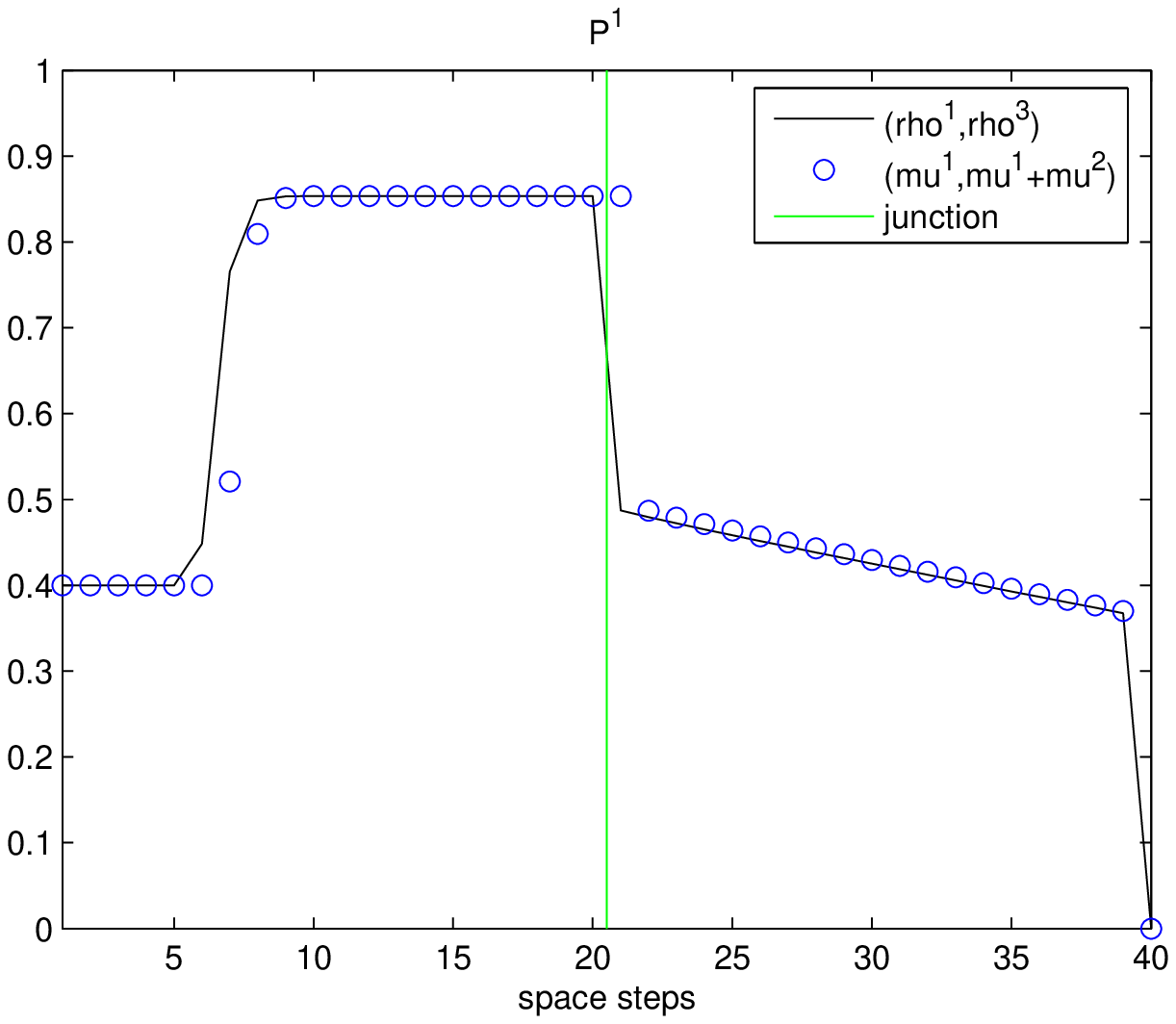}\\
\includegraphics[width=0.6\textwidth]{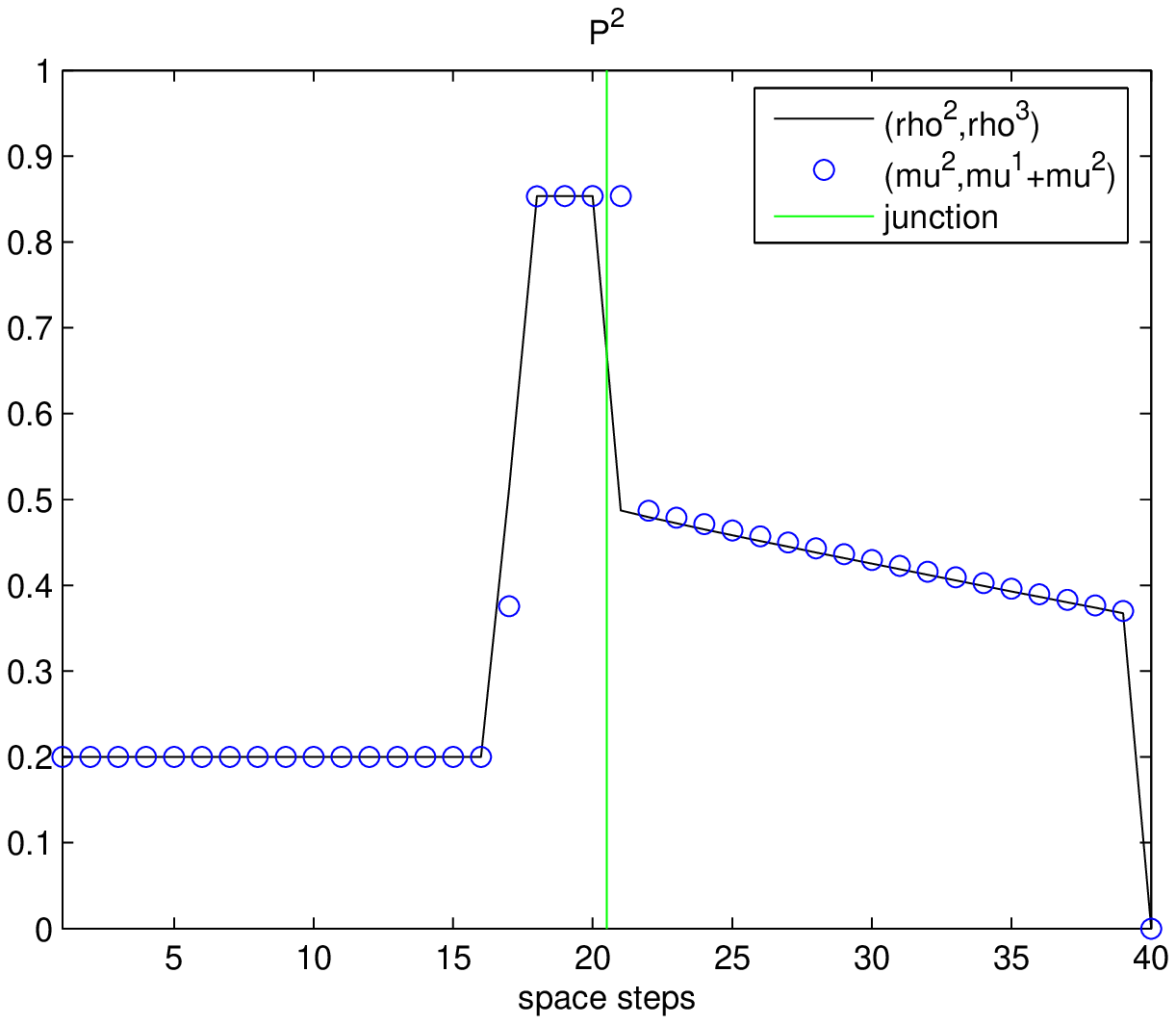}
\end{center}
\caption{Two incoming roads and one outgoing road. Numerical solution of classical and proposed algorithm with constant boundary conditions \eqref{2in1out_DBconst}. Path $P^1$ (top) and path $P^2$ (bottom). Note the 1-cell shift at the junction.}
\label{fig:2in1.confr}
\end{figure}

In the following we present the results of some numerical tests in which we compare the proposed algorithm \eqref{def:omega_approx}-\eqref{schema} with the classical one \eqref{Godunov}-\eqref{godunov_ben_in}. A theoretical comparison can be found in \cite{briani2013inprep}. For all tests we choose $\rhomax=1$ and $f(\rho)=\rho(1-\rho)$. We also assume that each arc has the same length, equal to 1, then the densities $\rho^1,\rho^2,\rho^3$ computed by the classical algorithm are defined in $[0,1]\times[0,t_f]$, while the densities $\mu^1,\mu^2$ computed by the proposed algorithm are defined in $[0,2]\times[0,t_f]$. We discretize each arc by means of 20 cells, then the junction is between the node $J-1=20$ and the node $J=21$. We choose $t_f=5$ and we discretize the time interval with 241 nodes. At time $t=0$ the network is empty. The boundary conditions at any time $t$ are:
\begin{equation}\label{2in1out_DBconst}
\begin{array}{l}
\rho^1(0,t)=\mu^1(0^{(1)},t)=0.4, \\
\rho^2(0,t)=\mu^2(0^{(2)},t)=0.2, \\
\rho^3(1,t)=\mu^1(2^{(1)},t)+\mu^2(2^{(2)},t)=0.
\end{array}
\end{equation}
By this choice, two queues are formed behind the junction, since the outgoing road is not able to receive the car flux coming from the two incoming roads.
Hereafter, we shall denote by $(\rho^i,\rho^j)$ the density on the path defined by two consecutive arcs $I_i$, $I_j$.
In Fig.\ \ref{fig:2in1.confr}
we report the solution at the final time computed by the classical algorithm with \textit{equidistributed priority coefficients} $q_1=q_2=1/2$, and that computed by the proposed algorithm.
For path $P^1$, the solution $(\rho^1,\rho^3)$ of the classical algorithm must be compared with the solution $(\mu^1,\mu^1+\mu^2)$ of the proposed algorithm, while for path $P^2$, the solution $(\rho^2,\rho^3)$ of the classical algorithm must be compared with the solution $(\mu^2,\mu^1+\mu^2)$ of the proposed algorithm. It can be seen that the two solutions are very similar, except for the fact that the junction is shifted by one cell. Indeed, \textit{the new algorithm makes the cell $J=21$ play the role of the last one of the incoming roads even if it is defined as the first cell of the outgoing road}. We also note that the proposed algorithm propagates the queue back in space slightly slower than the classical algorithm, probably because of the different behavior at the junction (see next test).

\begin{figure}[h!]
\begin{center}
\includegraphics[width=0.49\textwidth]{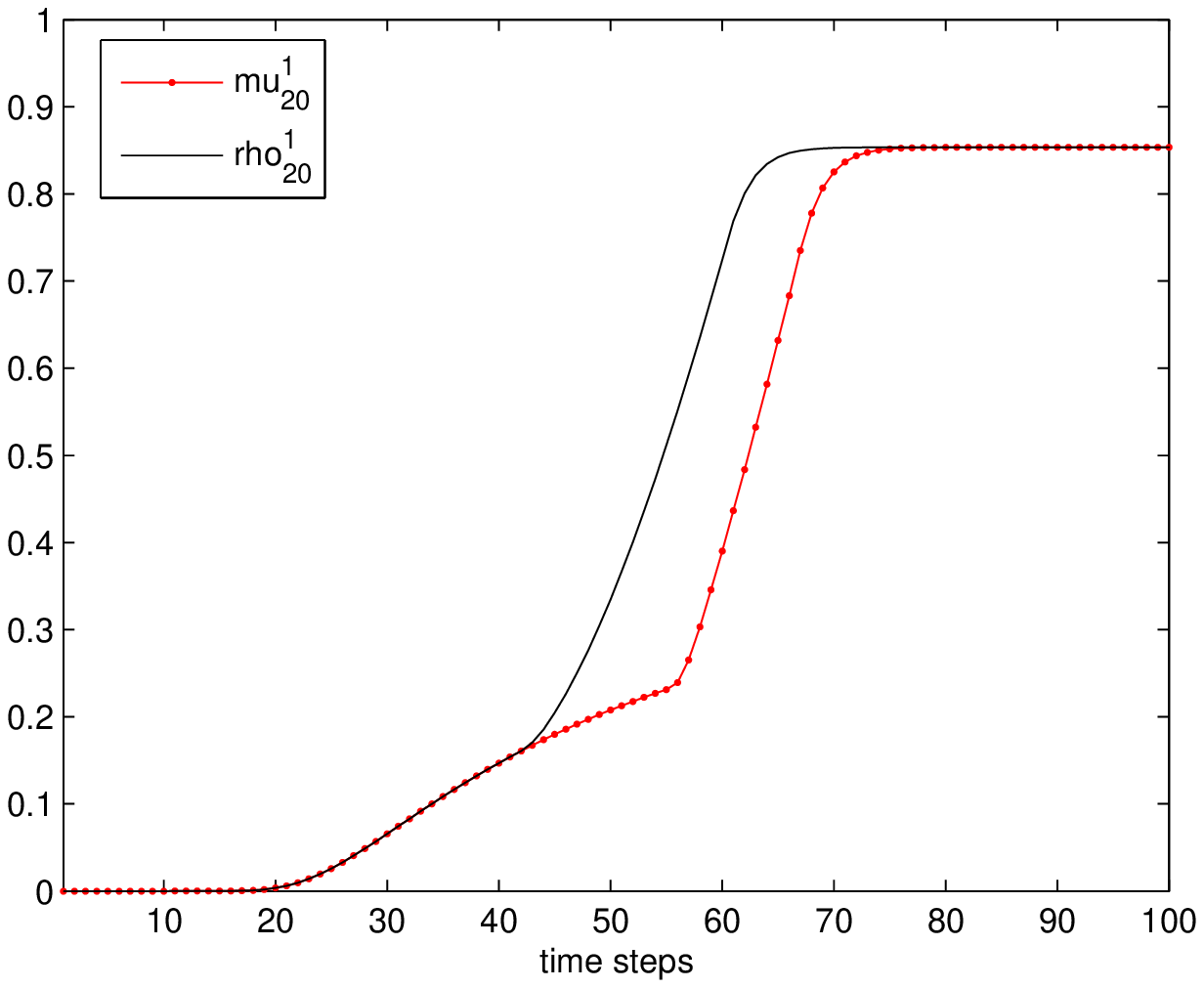}
\includegraphics[width=0.49\textwidth]{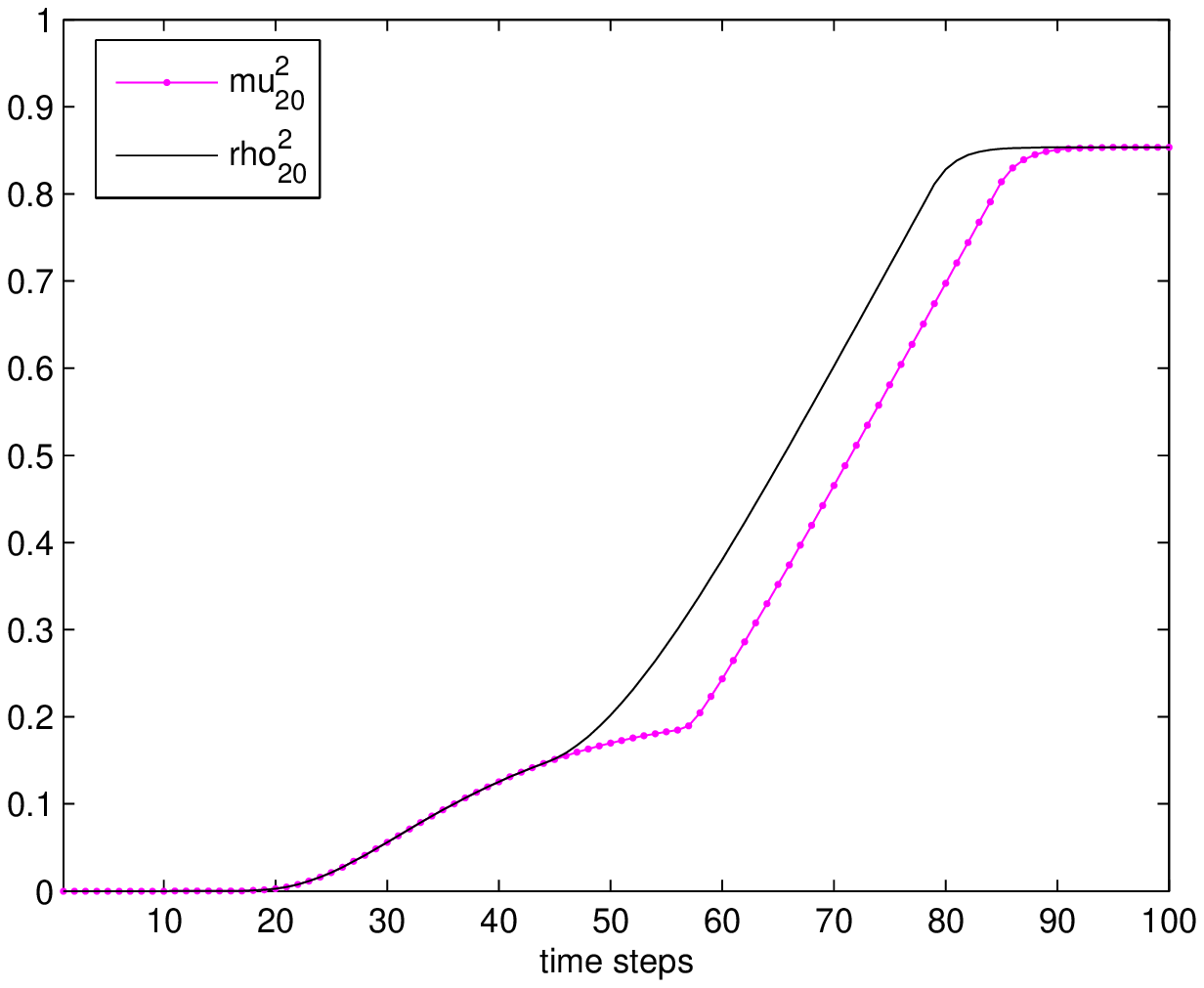}\\
\includegraphics[width=0.49\textwidth]{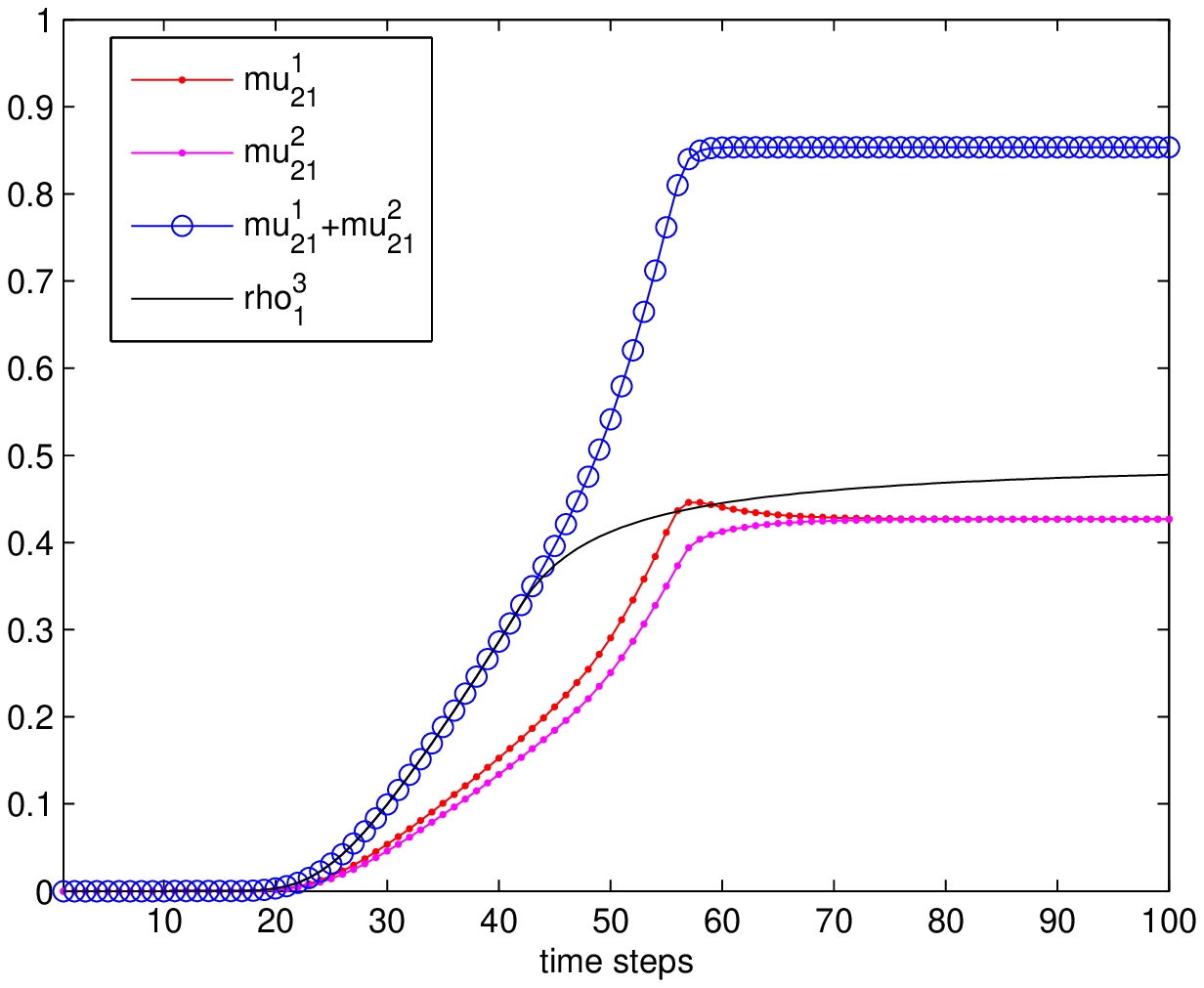}
\includegraphics[width=0.49\textwidth]{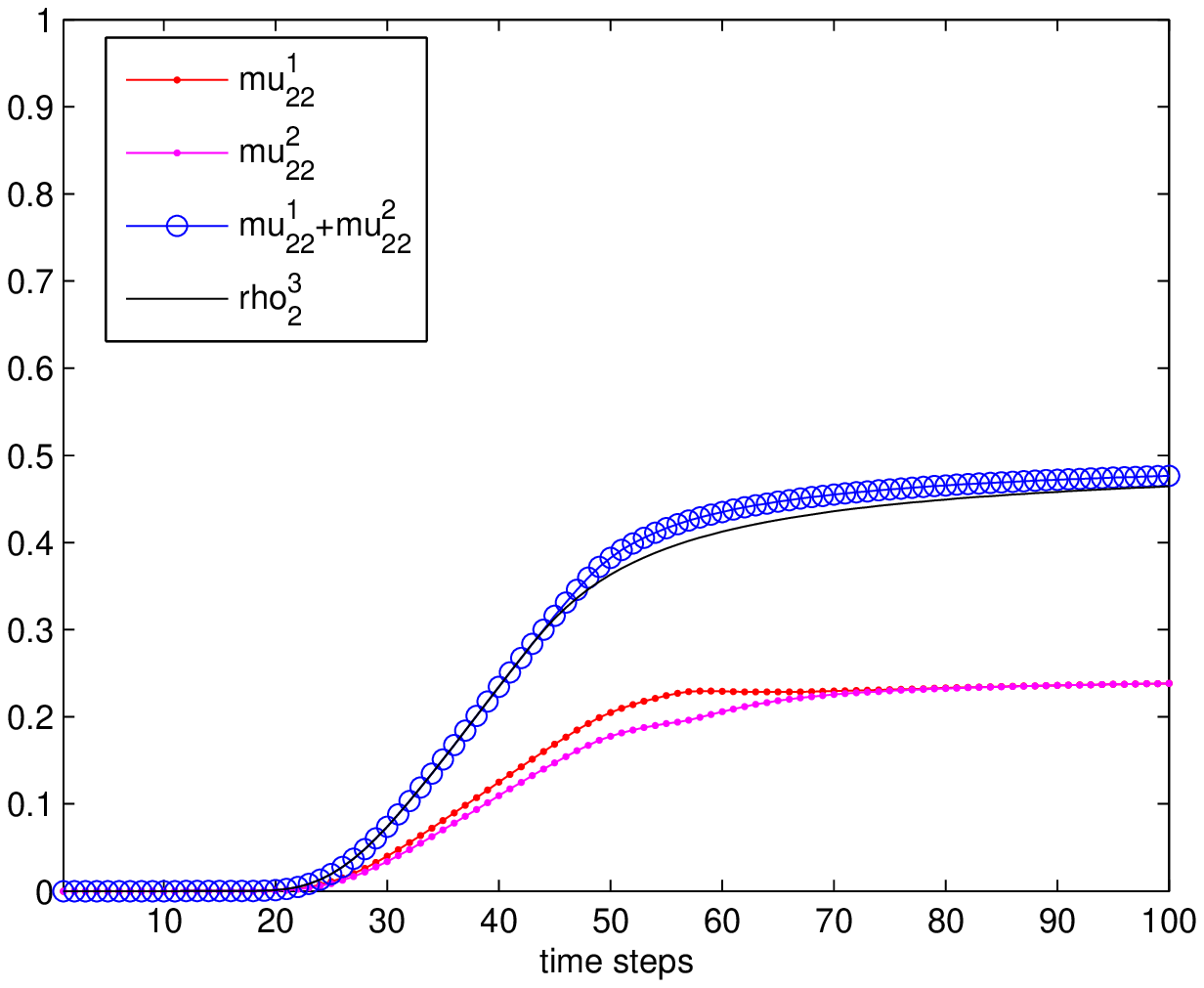}\\
\end{center}
\caption{Two incoming roads and one outgoing road. Evolution in time of the density around the junction, with boundary conditions \eqref{2in1out_DBconst}. Density on $P^1$ at cell $J-1=20$ (top-left), density on $P^2$ at cell $J-1=20$ (top-right), density at cell $J=21$ (bottom-left), and density at cell $J+1=22$ (bottom-right).}
\label{fig:2in1.trans}
\end{figure}

In Fig.\ \ref{fig:2in1.trans}
we show the evolution in time of the density at cells $k=J-1,J,J+1$. We note that, before the junction ($k=J-1$), the two algorithms do not coincide in a transient, in particular when the queue is forming and the characteristic curves change direction (from left$\rightarrow$right to right$\rightarrow$left). The proposed algorithm leads to a smaller density than the classical algorithm. A cell after the junction ($k=J+1$), instead, the situation is inverted: starting from the formation of the queue, the density of the proposed algorithm is slightly higher. Nevertheless, these differences are in some sense balanced, since the total number of vehicles which leave the network is exactly the same for the two algorithms, i.e.\
\begin{equation}\label{cfrINeOUT}
\int_0^{t_f}f(\rho^3(1-\Dx,t))dt=
\int_0^{t_f}f\big(\mu^1(2^{(1)}-\Dx,t)+\mu^2(2^{(2)}-\Dx)\big)dt.
\end{equation}
Just after the junction ($k=J$), the two algorithms show a remarkably different behavior, with two ``alien'' values $\mu^1_{J=21}$ and $\mu^2_{J=21}$ which are not immediately justified in the framework of the classical theory, see \cite{briani2013inprep}. These values cannot be considered by their own because they leave after the junction (where the two densities coexist). Rather, \textit{the sum of the two alien values is the value by means of which the two incoming roads perceive the presence of the junction ahead}.
The cell $J$ may be seen as a region with an oversize capacity which gathers the flows coming from the incoming roads. We may therefore say that the cell $J$ has been augmented by a ``buffer''.
However, our approach differs from the traffic models with buffer \cite{garavello2012DCDS, garavello2013bookchapt, herty2009NHM}. In fact, in those models the load of the buffer at any time is described by a dedicated function which evolves according to an additional ordinary differential equation.

\begin{figure}[h!]
\begin{center}
\includegraphics[width=0.49\textwidth]{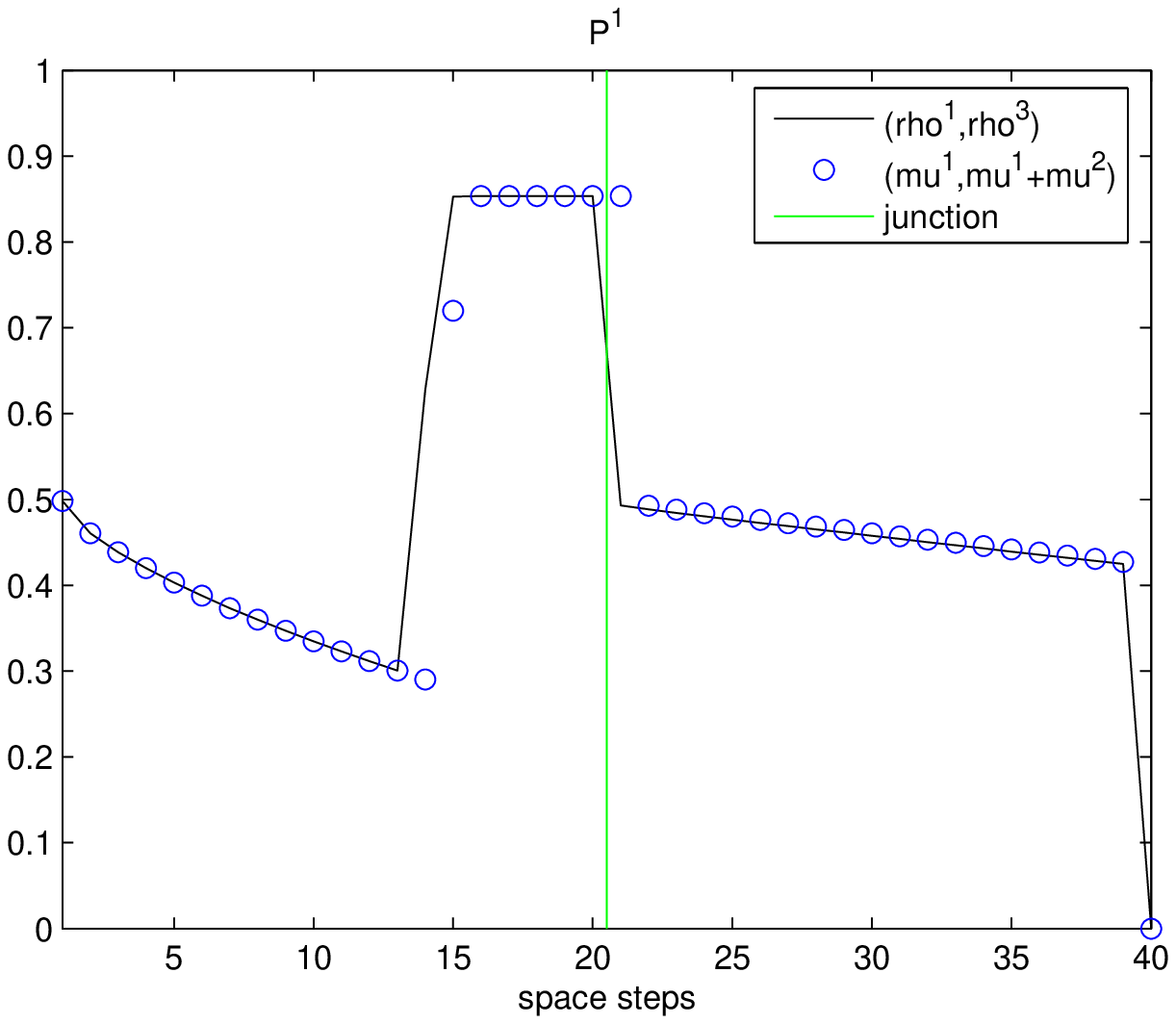}
\includegraphics[width=0.49\textwidth]{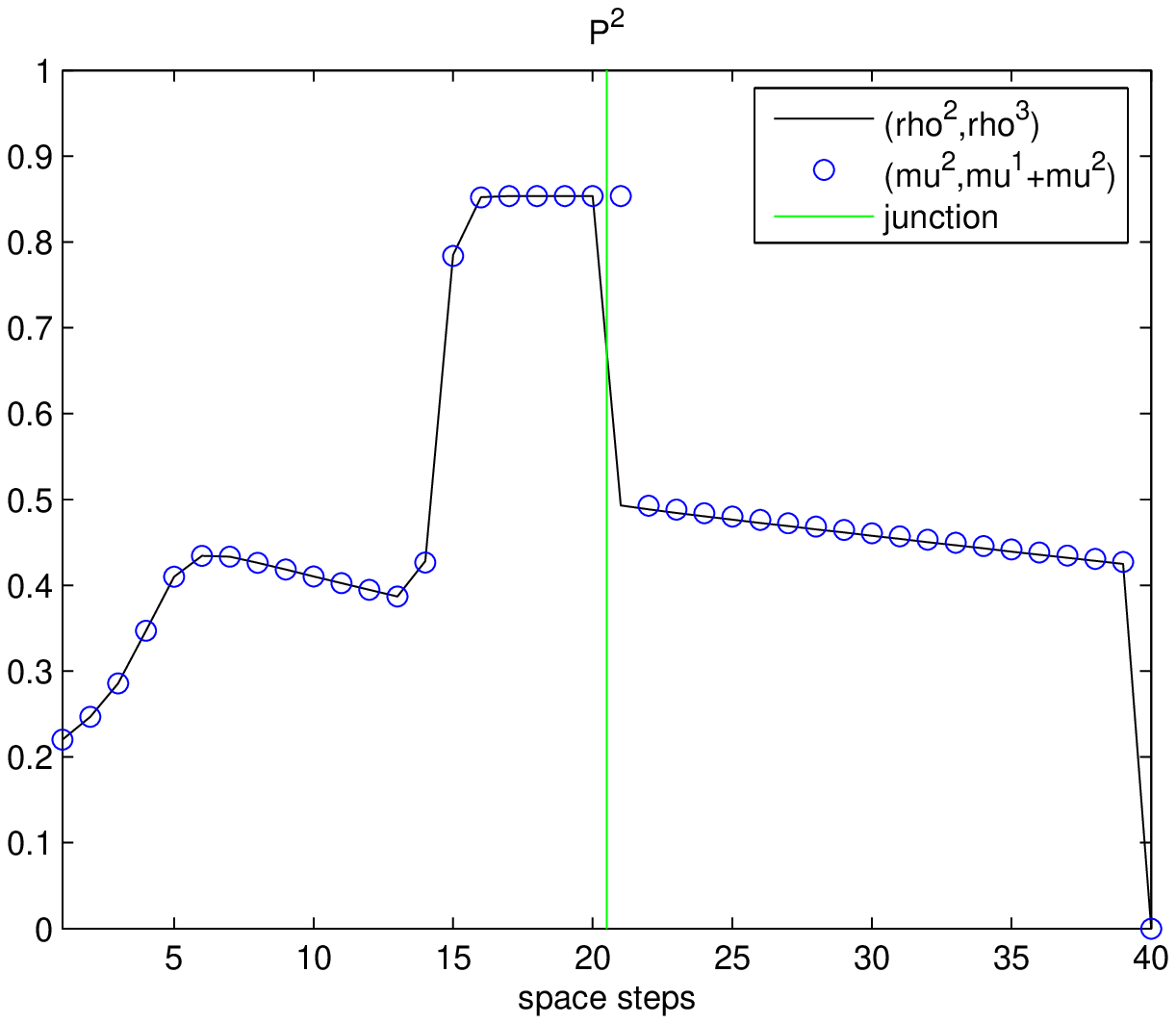}
\end{center}
\caption{Two incoming roads and one outgoing road. Numerical solution of classical and proposed algorithm with time-dependent left boundary conditions \eqref{2in1out_DBnonconst}. Path $P^1$ (left) and path $P^2$ (right). We observe the same behavior as in the previous test.}
\label{fig:2in1.DBnoncost}
\end{figure}

Finally, in Fig.\ \ref{fig:2in1.DBnoncost}
we report the solution at the final time $t_f=8$ computed by the classical algorithm with equidistributed priority coefficients $q_1=q_2=1/2$ and by the proposed algorithm in the case of a time-dependent left boundary conditions
\begin{equation}\label{2in1out_DBnonconst}
\begin{array}{l}
\rho^1(0,t)=\mu^1(0^{(1)},t)=\frac14(1+\sin t), \\
\rho^2(0,t)=\mu^2(0^{(2)},t)=\frac14(1+\cos t).
\end{array}
\end{equation}
In this test the transient period is kept alive by the time-dependent boundary conditions. Nevertheless, we observe the same behavior as in the previous test, meaning that the different solution around the junction do not affect the solution elsewhere.


\section{One incoming road and two outgoing roads}\label{sec:1in2out}
\begin{figure}[h!]
\begin{center}
\begin{psfrags}
\psfrag{P1}{$P^1$} \psfrag{P2}{$P^2$}
\psfrag{J}{\tiny $J$} \psfrag{Jp}{\tiny $J$\!+\!1} \psfrag{Jm}{\tiny $J$-1}
\psfrag{I1}{$[a_1,b_1]$} \psfrag{I2}{$[a_2,b_2]$} \psfrag{I3}{$[a_3,b_3]$}
\includegraphics[width=0.45\textwidth]{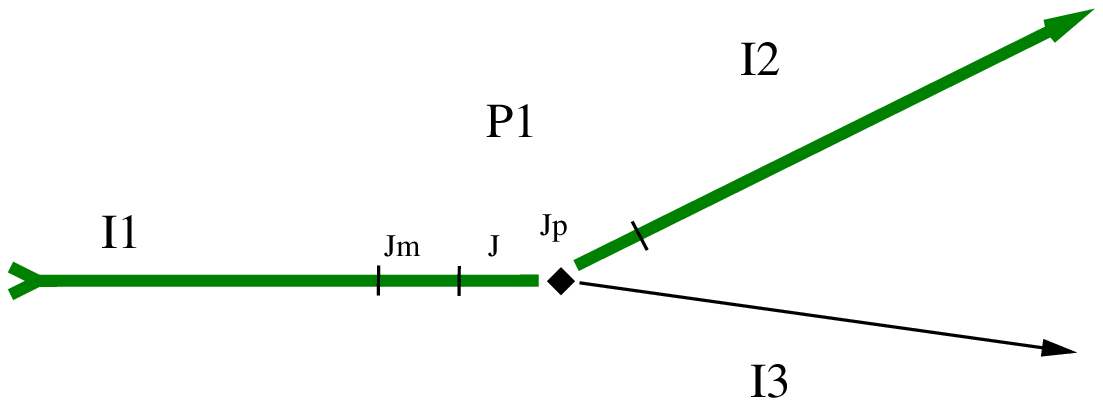}
\includegraphics[width=0.45\textwidth]{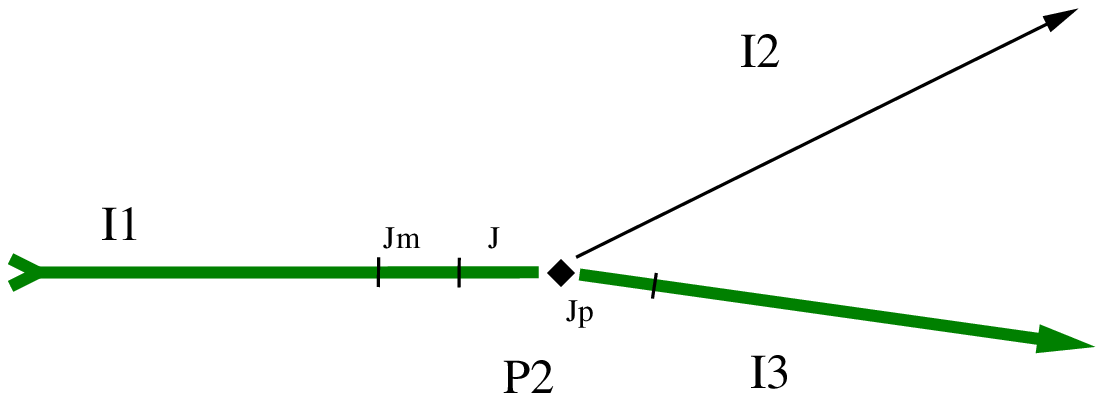}
\end{psfrags}
\end{center}
\caption{A network with 3 arcs and 1 junction. Path $P^1$ (left) and path $P^2$ (right).}
\label{fig:esempio_network_1in2out}
\end{figure}

In this section we consider a network with three arcs and one junction, with one incoming road and two outgoing roads. On the network two paths $P^1$ and $P^2$ are defined, see Fig.\ \ref{fig:esempio_network_1in2out}.
As before, we assume that each arc has the same length, equal to 1. We denote by $J$ the node \emph{just before} the junction, see Fig.\ \ref{fig:esempio_network_1in2out}. Note that the nodes after the junction, namely $J+1$, $J+2$, etc., can refer to one path or the other one, depending on the context.
We have
$$
\omega_{k^1}^{n}=\left\{\begin{array}{ll} \mu_{k^1}^{n,1}+\mu_{k^1}^{n,2} & k^1\leq J \\
\mu_{k^1}^{n,1} & k^1> J \end{array}\right.,
\qquad\qquad
\omega_{k^2}^{n}=\left\{\begin{array}{ll} \mu_{k^2}^{n,1}+\mu_{k^2}^{n,2} & k^2\leq J \\ 
\mu_{k^2}^{n,2} & k^2>J \end{array}\right.
$$
and the scheme (\ref{schema}) becomes \eqref{schema_2in1out}.

\begin{figure}[h!]
\begin{center}
\includegraphics[width=0.49\textwidth]{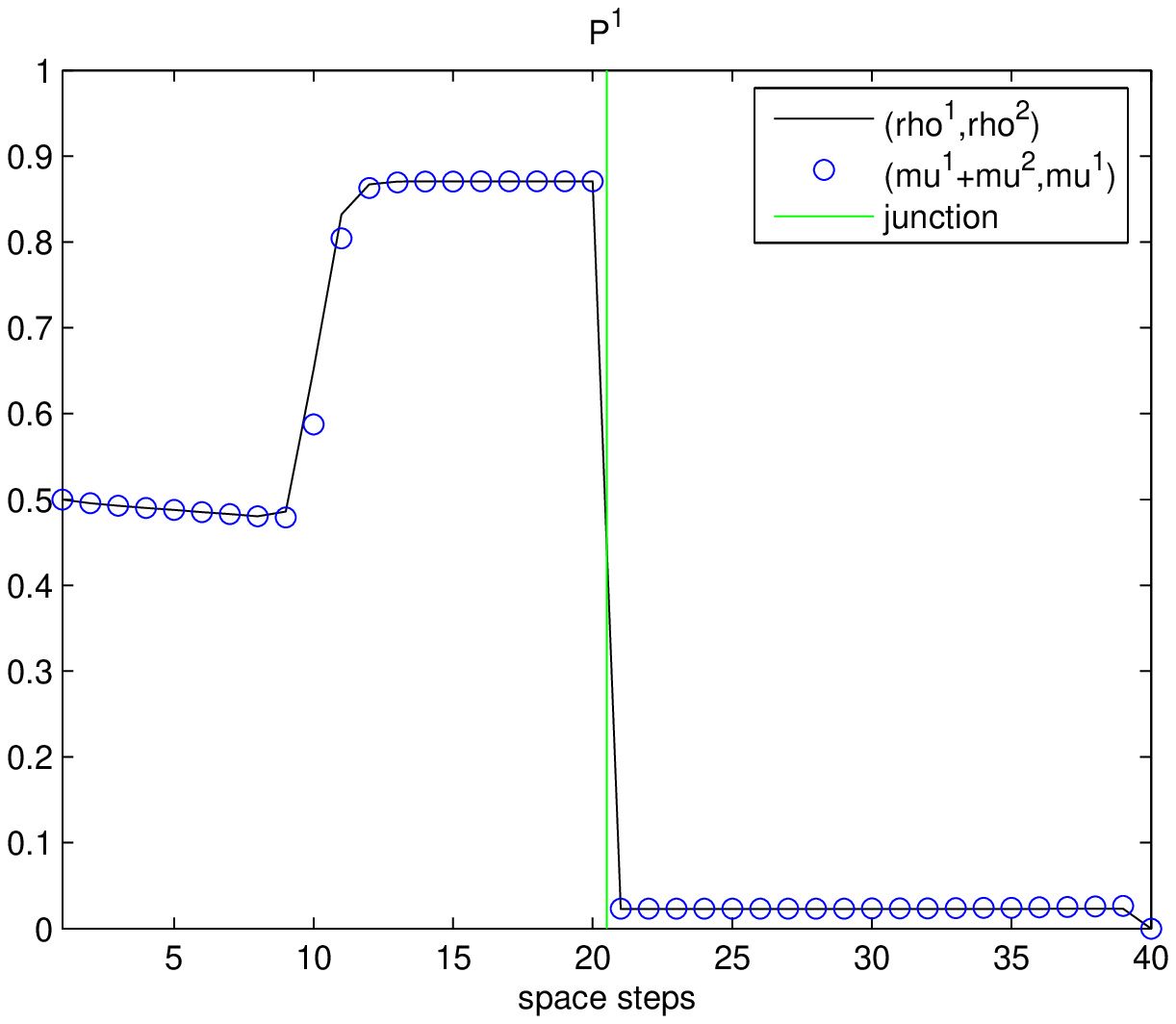}
\includegraphics[width=0.49\textwidth]{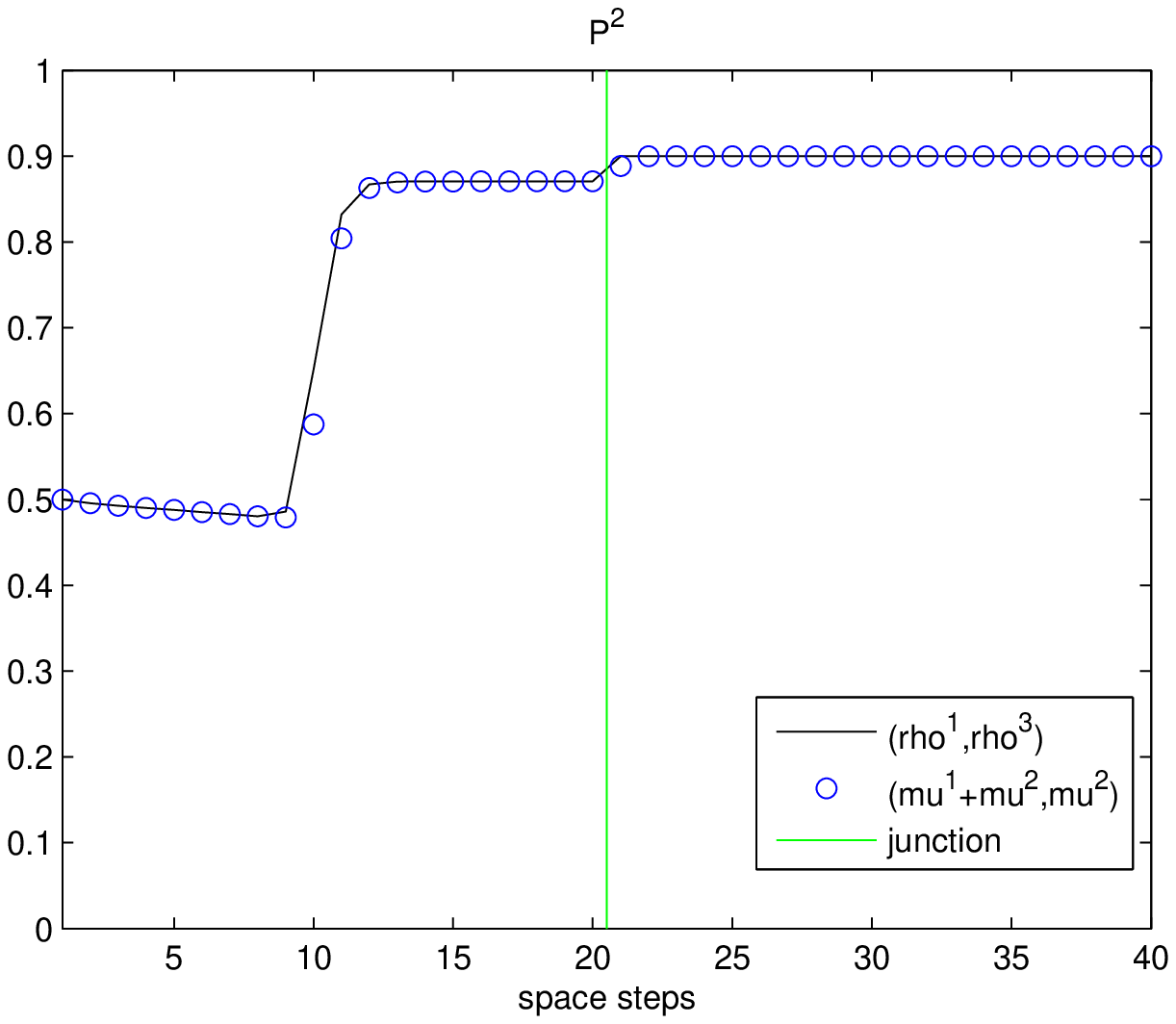}
\end{center}
\caption{One incoming road and two outgoing roads. Numerical solution of classical and proposed algorithm with boundary conditions (\ref{1in2out_DB_rho}),(\ref{1in2out_DB_mu}). Path $P^1$ (left) and path $P^2$ (right). The two solutions overlap.}
\label{fig:1in2.confr}
\end{figure}

As in the previous case, we discretize each arc by means of 20 cells. The junction is between the node $J=20$ and the node $J+1=21$. At time $t=0$ the network is empty. The boundary conditions for the classical algorithm are
\begin{equation}\label{1in2out_DB_rho}
\begin{array}{l}
\rho^1(0,t)=0.5, \quad \rho^2(1,t)=0, \quad \rho^3(1,t)=0.9
\end{array}
\end{equation}
and we set the preference coefficients (\ref{def:A}) as $\alpha_{21}=0.8$ and $\alpha_{31}=0.2$.
By this choice, a queue is formed behind the junction, since the road 3 is almost full. According to (\ref{1in2out_DB_rho}), we set
\begin{equation}\label{1in2out_DB_mu}
\begin{array}{l}
 \mu^1(0^{(1)},t)=\alpha_{21} \rho^1(0,t) = 0.4, \quad \mu^2(0^{(2)},t)=\alpha_{31} \rho^1(0,t) = 0.1, \\
 \mu^1(2^{(1)},t) = \rho^2(1,t), \quad \mu^2(2^{(2)},t)=\rho^3(1,t).
\end{array}
\end{equation}
In Fig.\ \ref{fig:1in2.confr}
we report the solution at the final time $t_f=11$ computed by the classical algorithm  and by the proposed algorithm. In this case the two solutions overlap, showing no shift in the localization of the junction.


\section{Two incoming roads and two outgoing roads}\label{sec:2in2out}
\begin{figure}[h!]
\begin{center}
\includegraphics[width=0.49\textwidth]{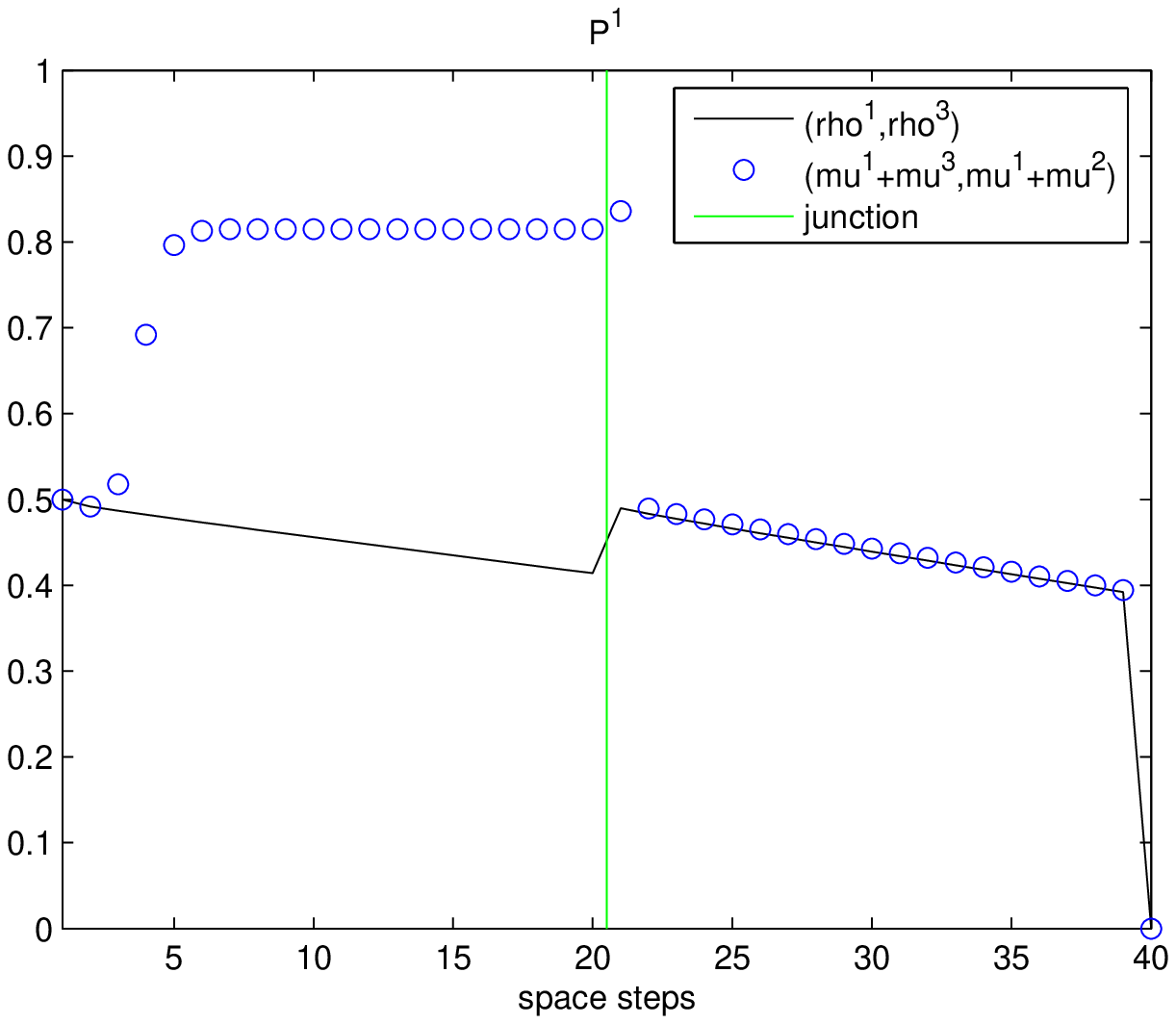}
\includegraphics[width=0.49\textwidth]{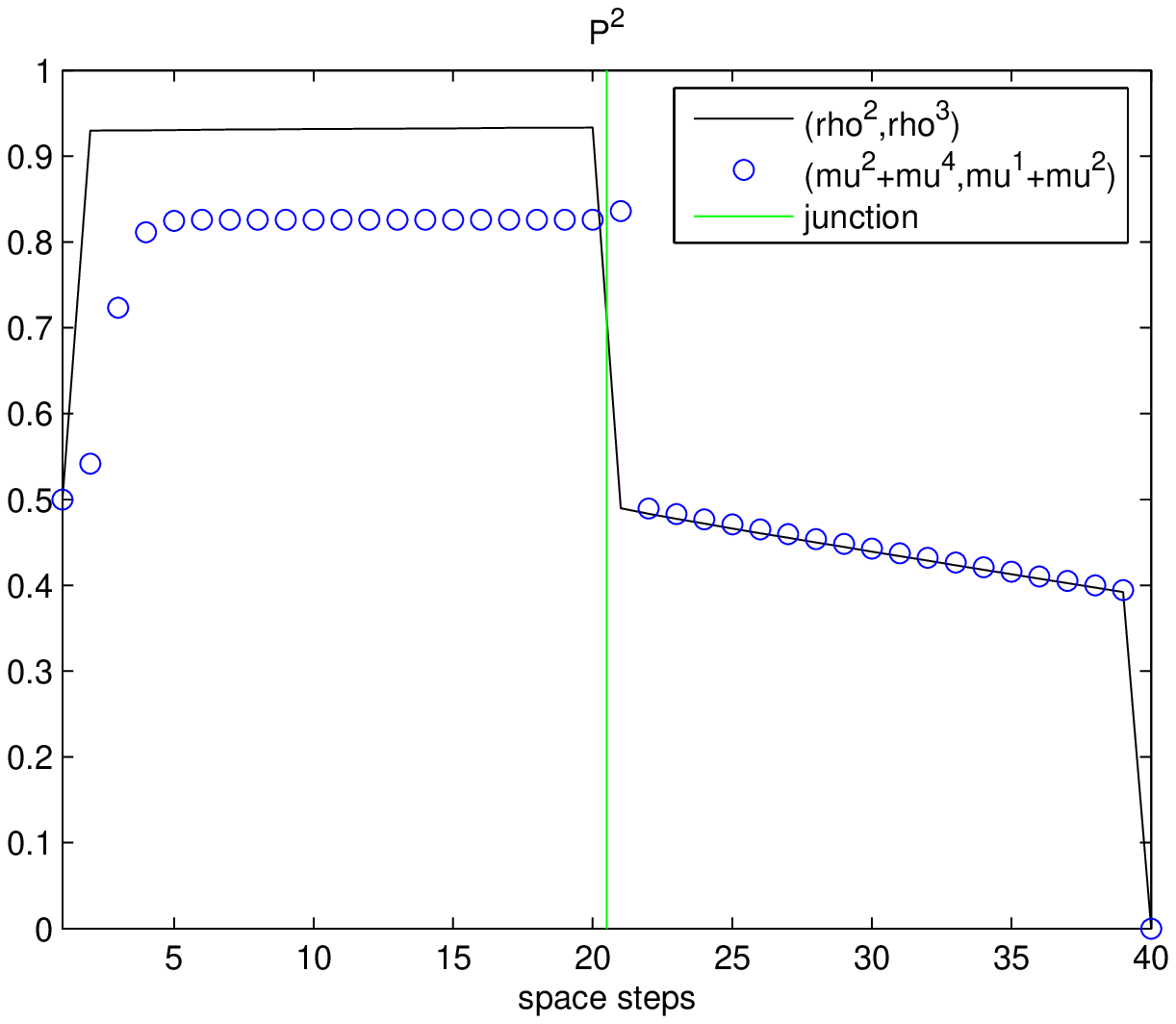}\\
\includegraphics[width=0.49\textwidth]{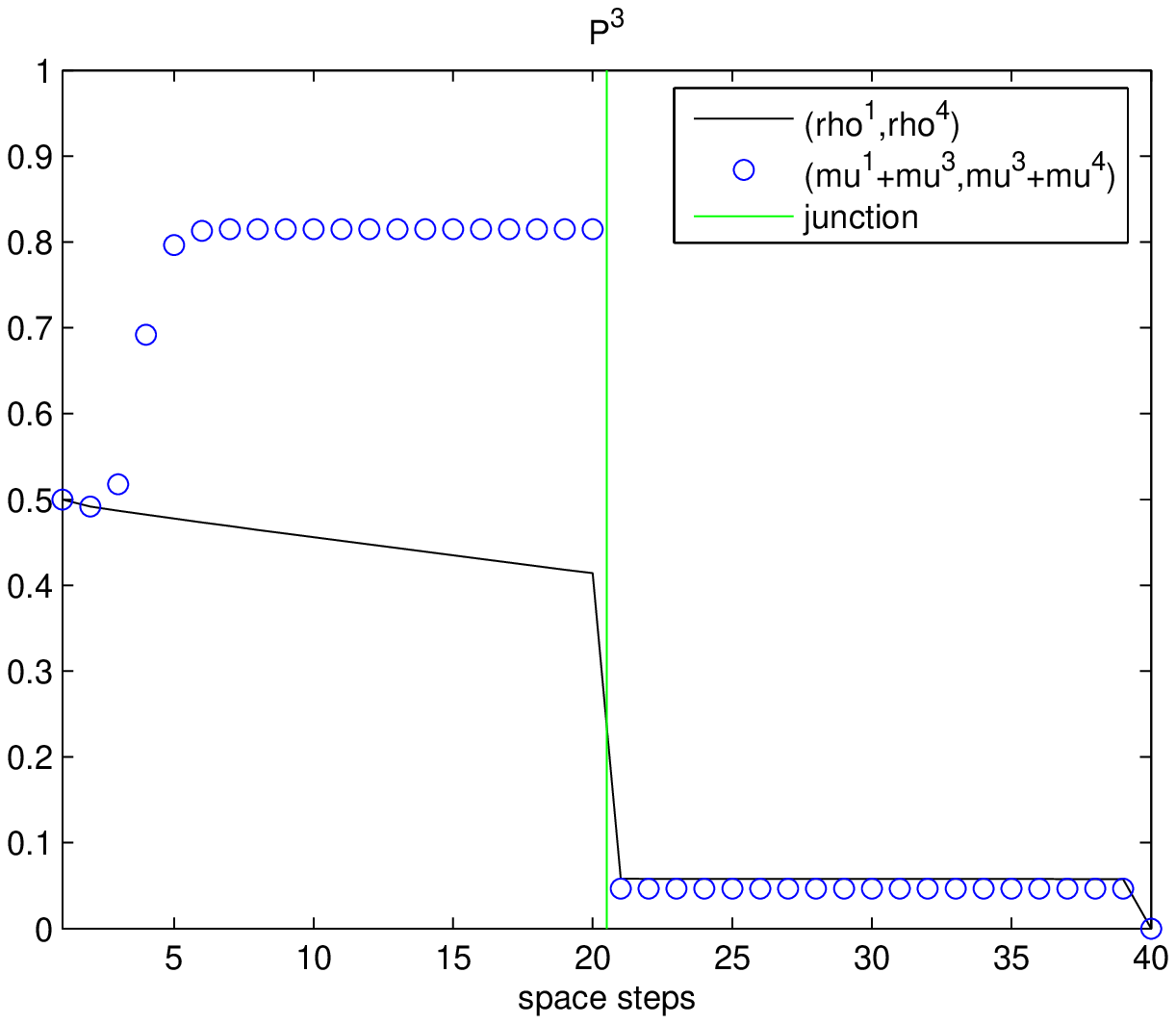}
\includegraphics[width=0.49\textwidth]{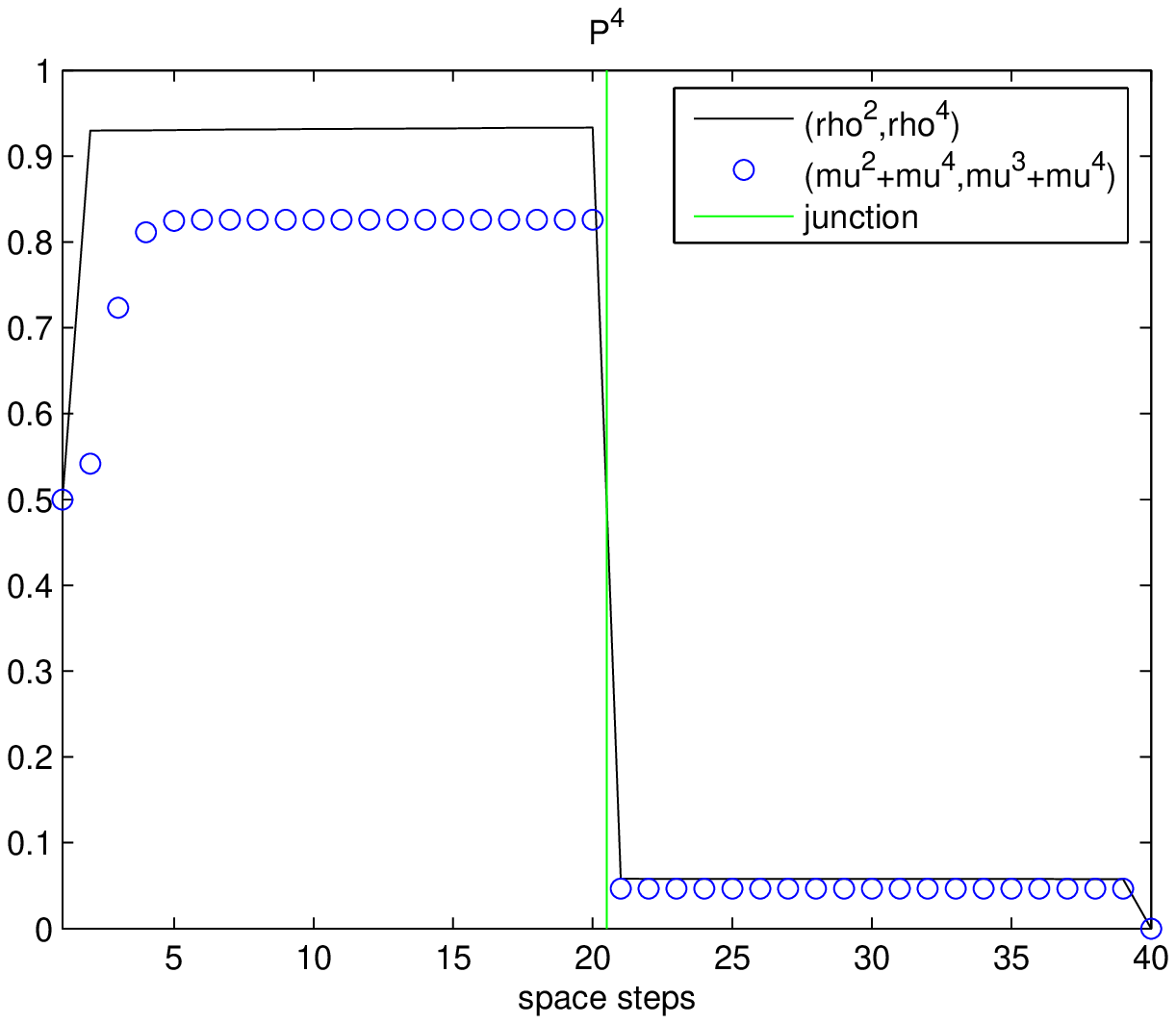}\\
\end{center}
\caption{Two incoming roads and two outgoing roads. Numerical solution of classical and proposed algorithm with boundary conditions \eqref{2in2out_DB}. Path $P^1$ (top-left), path $P^2$ (top-right), path $P^3$ (bottom-right), path $P^4$ (bottom-right). The two solutions are different.}
\label{fig:2in2}
\end{figure}

\begin{figure}[h!]
\begin{center}
\begin{psfrags}
\psfrag{g1}{$\gamma^1$} \psfrag{g2}{$\gamma^2$}
\psfrag{g1max}{$\gamma^1_{\textup{max}}$}
\psfrag{g2max}{$\gamma^2_{\textup{max}}$}
\psfrag{r1}{\tiny $\alpha_{41}\gamma^1\!+\!\alpha_{42}\gamma^2\!=\!\gamma^4_{\textup{max}}$}
\psfrag{r2}{\tiny $\alpha_{31}\gamma^1\!+\!\alpha_{32}\gamma^2\!=\!\gamma^3_{\textup{max}}$}
\psfrag{r3}{\tiny $\gamma^1\!-\!\gamma^2\!=\!0$}
\includegraphics[width=0.45\textwidth]{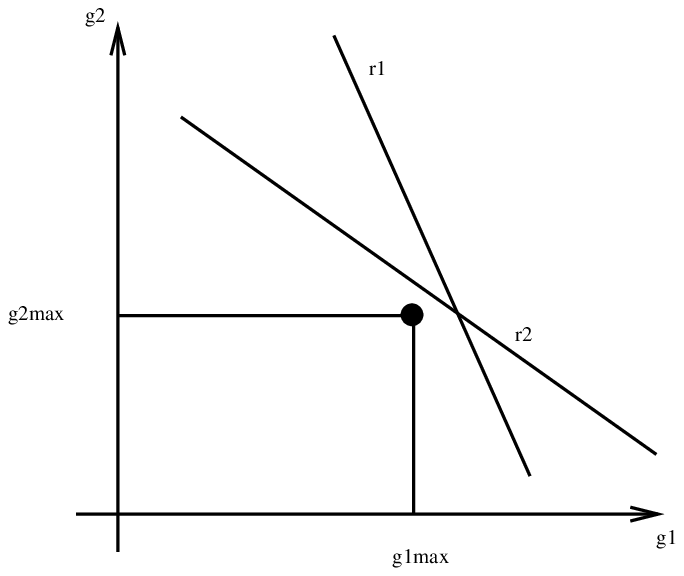} \qquad
\includegraphics[width=0.45\textwidth]{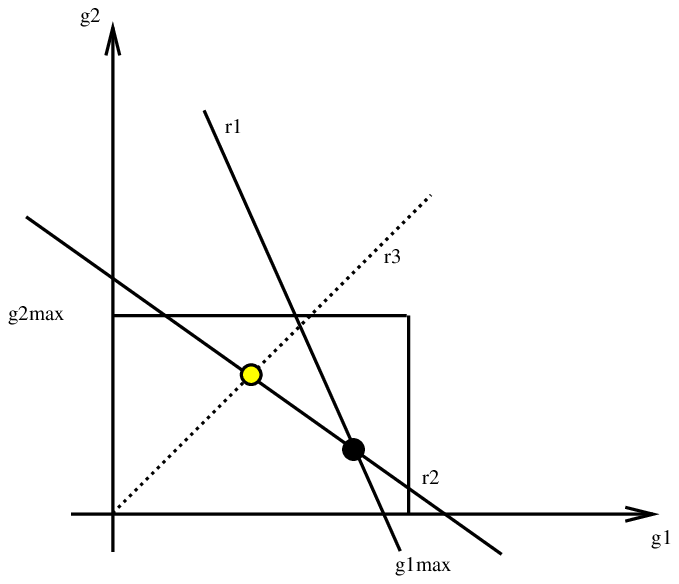}
\end{psfrags}
\end{center}
\caption{The solution of \eqref{PL} at junction found by maximizing the flux (black circle), and the one which seems to be found by the multi-path algorithm (yellow circle). In case of no queues the two solutions overlap (left), in case of queue formation the two solutions are different (right).}
\label{fig:2in2_PL}
\end{figure}

In this section we consider a network with four arcs and one junction, with two incoming roads and two outgoing roads. On the network four paths $P^1,P^2,P^3,P^4$ are defined. Boundary conditions and preference coefficients are:
\begin{equation}\label{2in2out_DB}
\begin{array}{l}
\rho^1(0,t)=\rho^2(0,t)=0.5, \qquad \rho^3(1,t)=\rho^4(1,t)=0, \\
\alpha_{31}=0.8, \quad \alpha_{41}=0.2, \quad \alpha_{32}=0.9, \quad \alpha_{42}=0.1, \\
 \mu^1(0^{(1)},t)=\alpha_{31}\rho^1(0,t), \quad \mu^2(0^{(2)},t)=\alpha_{32}\rho^2(0,t), \\ \mu^3(0^{(3)},t)=\alpha_{41}\rho^1(0,t), \quad \mu^4(0^{(4)},t)=\alpha_{42}\rho^2(0,t) , \\
  \mu^1(2^{(1)},t)=\mu^2(2^{(2)},t)=\mu^3(2^{(3)},t)=\mu^4(2^{(4)},t)=0.
\end{array}
\end{equation}
In Fig.\ \ref{fig:2in2}
we report the solution at the final time $t_f=6$ computed by the classical algorithm and by the proposed algorithm. In this case the solutions do not coincide. Also, the total number of vehicles that leave the network is not the same (cfr.\ \eqref{cfrINeOUT}), then we infer that the flux at junction must be different. As a consequence, \textit{the proposed algorithm does not maximize the flux at junction}. Remarkably, the proposed algorithm gives the same density in the incoming roads, independently from the choice of the boundary conditions. The question arises which solution is chosen by the proposed algorithm. Numerical evidence shows that the multi-path algorithm coincides with the classical one in case of no formation of queues. If, instead, the constraint $A(\gamma^1,\gamma^2)^T\in(\Omega_{3}\times\Omega_{4})$ in (\ref{Omega}) is active, the multi-path algorithm seems to find the solution of the problem (\ref{PL}), but with the additional (and not required for uniqueness of the solution) constraint (\ref{constraint_priorities}), assuming $q_1=q_2$, see Fig.\ \ref{fig:2in2_PL}.
%


\section{A real application}\label{sec:realapplication}
\begin{figure}[h!]
\begin{center}
\includegraphics[width=0.9\textwidth]{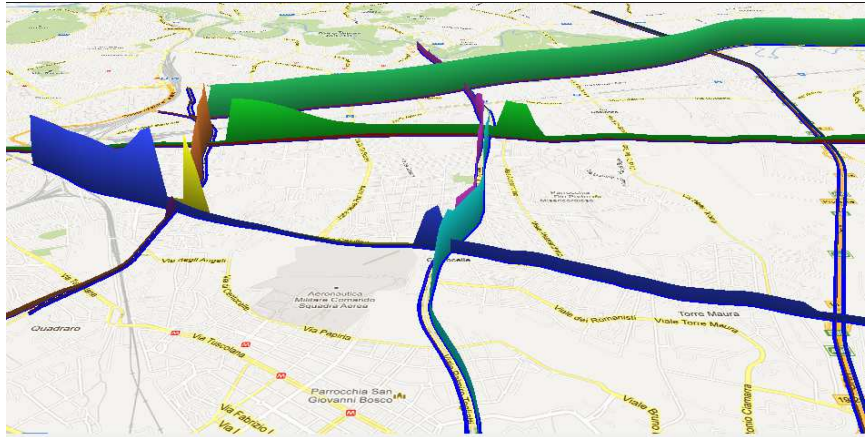}
\end{center}
\caption{A screenshot of the simulator. The density is represented as ``walls'' along the roads.}
\label{fig:ris_zeropiu}
\end{figure}

In order to test the proposed algorithm on a real case, we considered a network located in Rome (Italy), constituted by 6 two-lane roads and 7 junctions.
The whole network is 328.2 km.
%
%
We used the local version of the model described in Section \ref{sec:localmodel} and we set $\Delta x=100$ m, $\Delta t=2.5$ s. We also placed four traffic lights coordinated in pairs. The final time for the simulation is $t_f=\frac34$ h, see Fig.\ \ref{fig:ris_zeropiu} for a screenshot.
The code is written in C++ (serial) and run on an Intel i3 2.27 GHz processor. The CPU time for the entire simulation was 0.5 s. This result suggests that the proposed technique can be actually used to forecast traffic flow in large networks, keeping to a minimum the implementing effort.

\section{Conclusions and future work}
The outcomes of the numerical tests allow us to sketch some preliminary considerations. Avoiding any additional procedures at junctions (i.e.\ the most annoying part of existing algorithms), the multi-path algorithm is able to compute a solution which is admissible in the framework of the classical theory but does not provide in general the maximization of the flux at junctions. It is useful to stress here that the maximization of flux is only one of the possible conditions to ensure uniqueness of the solution. Moreover, in some cases this condition may prescribe strongly unbalanced fluxes from the incoming roads. The multi-path algorithm, instead, overrides the flux maximization favouring the flux balance.

We also stress that, once the topology of the network is suitably implemented in the function $\omega$ (total density), adding or removing a single road or path requires minimal effort.


\section*{Acknowledgments}
The authors thank R. Natalini and B. Piccoli for their useful suggestions, and R. M. Colombo for having pointed out the reference \cite{mercier2009JMAA}.

\medskip
\medskip


\begin{thebibliography}{99}

\bibitem{adreianov2011ARMA} (MR2807133) 
\newblock B. Andreianov, K. H. Karlsen and N. H. Risebro,
\newblock \emph{\emph{A theory of $L^1$-dissipative solvers for scalar conservation laws with discontinuous flux}},
\newblock \emph{Arch. Rational Mech. Anal.}, \textbf{201} (2011), 27--86.

\bibitem{benzoni2003EJAM} (MR2020123) 
\newblock S. Benzoni-Gavage and R. M. Colombo,
\newblock \emph{\emph{An $n$-populations model for traffic flow}},
\newblock \emph{European J. Appl. Math.}, \textbf{14} (2003), 587--612.

\bibitem{bressanbook} (MR1816648)
\newblock A. Bressan,
\newblock \emph{Hyperbolic Systems of Conservation Laws. The One-Dimensional Cauchy Problem},
\newblock Oxford Lecture Series in Mathematics, 20, Oxford University Press, New York, 2000.

\bibitem{bretti2007ACME} (MR2339915) 
\newblock G. Bretti, R. Natalini and B. Piccoli,
\newblock \emph{\emph{A fluid-dynamic traffic model on road networks}},
\newblock \emph{Arch. Comput. Methods Eng.}, \textbf{14} (2007), 139--172.

\bibitem{briani2013inprep}
\newblock M. Briani and E. Cristiani,
\newblock \emph{\emph{An easy-to-use algorithm for simulating traffic flow on networks: theoretical study}},
\newblock Submitted, arXiv:1401.1651.

\bibitem{burger2008JEM} (MR2396483) 
\newblock R. B\"urger and K. H. Karlsen,
\newblock \emph{\emph{Conservation laws with discontinuous flux: A short introduction}},
\newblock \emph{J. Eng. Math.}, \textbf{60} (2008), 241--247.

\bibitem{chitour2005DCDS-B} (MR2151724) 
\newblock Y. Chitour and B. Piccoli,
\newblock \emph{\emph{Traffic circles and timing of traffic lights for cars flow}},
\newblock \emph{Discrete and Continuous Dynamical Systems - Series B}, \textbf{5} (2005), 599--630.

\bibitem{coclite2005SIMA} (MR2178224) 
\newblock G. M. Coclite, M. Garavello and B. Piccoli,
\newblock \emph{\emph{Traffic flow on a road network}},
\newblock \emph{SIAM J. Math. Anal.}, \textbf{36} (2005), 1862--1886.

\bibitem{cristiani2010CAIM} (MR2812252) 
\newblock E. Cristiani, C. de Fabritiis and B. Piccoli,
\newblock \emph{\emph{A fluid dynamic approach for traffic forecast from mobile sensor data}},
\newblock \emph{Commun. Appl. Ind. Math.}, \textbf{1} (2010), 54--71.

\bibitem{daganzo1994TRB} 
\newblock C. F. Daganzo,
\newblock \emph{\emph{The cell transmission model: A dynamic representation of highway traffic consistent with the hydrodynamic theory}},
\newblock \emph{Transportation Research Part B}, \textbf{28} (1994), 269--287.

\bibitem{daganzo1995TRBa} 
\newblock C. F. Daganzo,
\newblock \emph{\emph{The cell transmission model, part II: Network traffic}},
\newblock \emph{Transportation Research Part B}, \textbf{29} (1995), 79--93.

\bibitem{garavello2012DCDS} (MR2885791) 
\newblock M. Garavello and P. Goatin,
\newblock \emph{\emph{The Cauchy problem at a node with buffer}},
\newblock \emph{Discrete and Continuous Dynamical Systems - Series A}, \textbf{32} (2012), 1915--1938.

\bibitem{garavello2005CMS} (MR2165018) 
\newblock M. Garavello and B. Piccoli,
\newblock \emph{\emph{Source-destination flow on a road network}},
\newblock \emph{Comm. Math. Sci.}, \textbf{3} (2005), 261--283.

\bibitem{piccolibook} (MR2328174)
\newblock M. Garavello and B. Piccoli,
\newblock \emph{Traffic Flow on Networks},
\newblock AIMS Series on Applied Mathematics, 1, 2006.

\bibitem{garavello2013bookchapt} 
\newblock M. Garavello and B. Piccoli,
\newblock \emph{\emph{A multibuffer model for LWR road networks}},
\newblock \emph{Advances in Dynamic Network Modeling in Complex Transportation Systems, Complex Networks and Dynamic Systems}, \textbf{2} (2013), 143--161.

\bibitem{herrera2010TRB} 
\newblock J. C. Herrera and A. M. Bayen,
\newblock \emph{\emph{Incorporation of Lagrangian measurements in freeway traffic state estimation}},
\newblock \emph{Transportation Research Part B}, \textbf{44} (2010), 460--481.

\bibitem{herty2009NHM} (MR2552171) 
\newblock M. Herty, J.-P. Lebacque and S. Moutari,
\newblock \emph{\emph{A novel model for intersections of vehicular traffic flow}},
\newblock \emph{Netw. Heterog. Media}, \textbf{4} (2009), 813--826.

\bibitem{herty2007ANM} (MR2310753) 
\newblock M. Herty, M. Sea\"id and A. K. Singh,
\newblock \emph{\emph{A domain decomposition method for conservation laws with discontinuous flux function}},
\newblock \emph{Appl. Numer. Math.}, \textbf{57} (2007), 361--373.

\bibitem{hilliges1995TRB} 
\newblock M. Hilliges and W. Weidlich,
\newblock \emph{\emph{A phenomenological model for dynamic traffic flow in networks}},
\newblock \emph{Transportation Research Part B}, \textbf{29} (1995), 407--431.

\bibitem{holden1995SIMA} (MR1338371) 
\newblock H. Holden and N. H. Risebro,
\newblock \emph{\emph{A mathematical model of traffic flow on a network of unidirectional roads}},
\newblock \emph{SIAM J. Math. Anal.}, \textbf{26} (1995), 999--1017.

\bibitem{lebacque1996proc}
\newblock J.-P. Lebacque,
\newblock \emph{The Godunov scheme and what it means for first order traffic flow models},
\newblock in \emph{Proc. of the 13th International Symposium on Transportation and Traffic Theory, Lyon, France} (ed. J. B. Lesort), Elsevier, 1996, 647--677.

\bibitem{levequebook} (MR1153252) 
\newblock R. J. LeVeque,
\newblock \emph{Numerical Methods for Conservation Laws},
\newblock Birkh\"auser, 1992.

\bibitem{LW} (MR0072606) 
\newblock M. J. Lighthill and G. B. Whitham,
\newblock \emph{\emph{On kinetic waves. II. Theory of traffic flows on long crowded roads}},
\newblock \emph{Proc. Roy. Soc. Lond. A}, \textbf{229} (1955), 317--345.

\bibitem{mercier2009JMAA} (MR2476922) 
\newblock M. Mercier,
\newblock \emph{\emph{Traffic flow modelling with junctions}},
\newblock \emph{J. Math. Anal. Appl.}, \textbf{350} (2009), 369--383.

\bibitem{R} (MR0075522) 
\newblock P. I. Richards,
\newblock \emph{\emph{Shock waves on the highway}},
\newblock \emph{Operation Research}, \textbf{4} (1956), 42--51.

\bibitem{towers2000SINUM} (MR1770068) 
\newblock J. D. Towers,
\newblock \emph{\emph{Convergence of a difference scheme for conservation laws with a discontinuous flux}},
\newblock \emph{SIAM J. Numer. Anal.}, \textbf{38} (2000), 681--698.

\bibitem{wong2002TRA} 
\newblock G. C. K. Wong and S. C. Wong,
\newblock \emph{\emph{A multi-class traffic flow model -- an extension of LWR model with heterogeneous drivers}},
\newblock \emph{Transportation Research Part A}, \textbf{36} (2002), 827--841.



\end{thebibliography}
\end{document}